\documentclass[12pt,english]{article}
\usepackage{mathptmx}
\usepackage{geometry}
\usepackage{color}
\usepackage{array}
\usepackage{bbding}
\usepackage{multirow}
\usepackage[fleqn]{amsmath}
\usepackage{amssymb}
\usepackage{amsthm}
\usepackage{graphicx}
\usepackage{natbib}
\usepackage{lscape}
\usepackage[hyphens]{url}
\usepackage[hidelinks]{hyperref}
\usepackage{comment}
\usepackage{bm}
\usepackage[title]{appendix}
\usepackage{longtable}
\usepackage{mathtools}
\usepackage{chngcntr}
\usepackage{authblk}
\usepackage{babel}
\usepackage{dsfont}
\usepackage{subcaption}
\usepackage{fancyhdr}
\usepackage{cleveref}
\usepackage{textgreek}
\usepackage{algorithm,algorithmic}
\usepackage[breakable]{tcolorbox}

\makeatletter
\allowdisplaybreaks
\date{}

\makeatletter
\def\blfootnote{\xdef\@thefnmark{}\@footnotetext}
\makeatother

\usepackage[mathlines, pagewise]{lineno}
\usepackage{etoolbox} 

\newcommand*\linenomathpatchAMS[1]{%
	\expandafter\pretocmd\csname #1\endcsname {\linenomathAMS}{}{}%
	\expandafter\pretocmd\csname #1*\endcsname{\linenomathAMS}{}{}%
	\expandafter\apptocmd\csname end#1\endcsname {\endlinenomath}{}{}%
	\expandafter\apptocmd\csname end#1*\endcsname{\endlinenomath}{}{}%
}

\expandafter\ifx\linenomath\linenomathWithnumbers
\let\linenomathAMS\linenomathWithnumbers
\patchcmd\linenomathAMS{\advance\postdisplaypenalty\linenopenalty}{}{}{}
\else
\let\linenomathAMS\linenomathNonumbers
\fi

\linenomathpatchAMS{gather}
\linenomathpatchAMS{multline}
\linenomathpatchAMS{align}
\linenomathpatchAMS{alignat}
\linenomathpatchAMS{flalign}

\newtheorem{prop}{Proposition}

\crefname{equation}{Eq.}{Eqs.}

\AddToHook{cmd/appendix/before}{%
  \setcounter{equation}{0}%
  \renewcommand{\theequation}{A.\arabic{equation}}%
  \setcounter{figure}{0}%
  \renewcommand{\thefigure}{A.\arabic{figure}}%
\setcounter{table}{0}%
  \renewcommand{\thetable}{A.\arabic{table}}%
  \setcounter{theorem}{0}%
  \setcounter{prop}{0}%
  \setcounter{lemma}{0}%
  \setcounter{cor}{0}%
}

\newcommand{\netT}{\mathcal{T}}

\newcommand{\netX}{\mathcal{X}}
\newcommand{\netY}{\mathcal{Y}}
\newcommand{\netZ}{\mathcal{Z}}

\newcommand{\E}{\mathop{{}\mathbb{E}}}

\begin{document}
	
\title{Dynamic Order Fulfillment in Last-Mile Drone Delivery under Demand Uncertainty}
\author[1]{Linxuan Shi}
\author[1]{Zhengtian Xu\footnote{Corresponding author. E-mail address: \textcolor{blue}{zhengtian@gwu.edu} (Z. Xu).}}
\author[2]{Miguel A. Lejeune}
\author[3]{Qi Luo}
\affil[1]{\small\emph{Department of Civil and Environmental Engineering, School of Engineering and Applied Science, George Washington University, Washington, DC, United States}\normalsize}
\affil[2]{\small\emph{Department of Decision Sciences, School of Business, George Washington University, Washington, DC, United States}\normalsize}
\affil[3]{\small\emph{Department of Business Analytics, Tippie College of Business, University of Iowa, Iowa City, IA, United States}\normalsize}

\date{\today}
\maketitle

\vspace{-0.6cm} 

\begin{abstract}
\noindent Drones have attracted growing interest in last-mile delivery due to their potential to significantly reduce costs and enhance operational flexibility, particularly in areas of sparse and uncertain demand where traditional truck delivery proves inefficient. This paper addresses the dynamic order fulfillment problem faced by a retailer operating a fleet of drones to service delivery requests that arrive stochastically. These delivery requests may vary in package profiles, delivery locations, and urgency. We adopt a rolling-horizon framework for order fulfillment and devise a two-stage stochastic program aimed at strategically managing existing orders while considering incoming requests that are subject to various uncertainties. A significant challenge in deploying the envisioned two-stage model lies in its incorporation of vehicle routing constraints, on which exact or brute-force methods are computationally inefficient and unsuitable for real-time operational decisions. To address this, we propose an accelerated L-shaped algorithm that (i) reduces the branching tree size, (ii) replaces exact second-stage solutions with heuristic estimates, and (iii) adapts an alternating strategy for adding optimality cuts. The proposed heuristic demonstrates remarkable performance superiority over the exact method, achieving a 20-fold reduction in average runtime while maintaining an average optimality gap of less than 1\%. We apply the algorithm to a wide range of instances to evaluate the benefits of postponing orders for batch service using the stochastic model. Our results show potential long-term cost savings of up to 20\% when demand uncertainty is explicitly considered in order fulfillment decisions. Meanwhile, the derived savings tend to diminish as the uncertainty increases in order arrivals.
\end{abstract}


\newpage
\section{Introduction}\label{sec:introduction}

Last-mile delivery has become an increasingly critical component of the logistics and transportation industry. In recent years, the COVID-19 pandemic prompted more consumers to favor online shopping over traditional retail methods, including those in rural areas who were previously less engaged with online channels \citep{shi2024dine}. As a result, the e-commerce industry is poised to capture a 41\% share of global retail sales by 2027, a substantial increase from its 18\% share in 2017 \citep{barthel2023winning}, which will continue to expand the frontier of last-mile delivery and place significant strain on logistics companies \citep{galiullina2024demand}. While urban areas grapple with congestion and high delivery volumes, the challenges in rural and remote regions are markedly different. Low population density results in fewer delivery stops spread across vast geographic areas. The spatiotemporal sparsity of demand makes rural logistics inherently costly to operate, often requiring companies to keep vehicles on standby for sporadic orders. Traditional truck-based delivery models, which excel in dense areas by enabling economies of scale, become particularly inefficient in scenarios of low and uncertain demand. With few drop-offs over long-distance routes, trucks are underutilized, significantly increasing fuel, labor, and maintenance costs per package. Additionally, poor road conditions and extended travel times further exacerbate operational inefficiencies, making conventional truck delivery neither practical nor cost-effective in sparsely populated areas \citep{fulfillment2025overcome}.

By contrast, drones or Unmanned aerial vehicles (UAVs) offer lower per-mile cost \citep{carlsson2018coordinated, she2024hybrid}, reduce delivery time \citep{sudbury2016cost, jennings2019study}, and decrease fuel consumption and emissions \citep{figliozzi2020carbon} compared to traditional truck deliveries. These advantages enable drones to effectively address many of the logistical challenges associated with sparse and uncertain demand. Moreover, their high speed and on-demand capability make them particularly well-suited for sporadic deliveries, such as urgent medical supplies, without requiring full-scale vehicle dispatch. Owing to their unparalleled convenience for consumers, drone-based last-mile delivery has steadily gained momentum over the past decade, expanding to a wide range of products and services, including goods, meals, groceries, and medications \citep{liu2019optimization, yin2024exact}. A growing number of delivery companies and online retailers, including FedEx and Amazon, are joining the game and testing drones for their last-mile delivery solutions \citep{reuters2022amazon}. Notably, companies such as the e-commerce platform JD and the logistics provider SF Express have prioritized rural villages in China, where underdeveloped ground transportation infrastructure presents a natural opportunity for drones to bridge the accessibility gap \citep{TechnologyReview2023}. In these emerging service contexts, the shift toward an on-demand operational paradigm introduces a need for more real-time and time-sensitive decision-making in autonomous order fulfillment. Pre-disclosed deterministic arrival of requests and delivery locations, which is prevalently assumed by most existing studies in this area, may no longer be applicable. Therefore, dynamic processing of order fulfillment and enabling the efficient batching of orders arriving at different time points in uncertain market contexts become crucial in the evolving landscape of last-mile delivery. To aid such operations, this study solves the dynamic order fulfillment problem for on-demand last-mile delivery using drones.

Consider the typical retail logistics where a service depot, such as a restaurant or pharmacy, continuously receives delivery orders, each distinguished by varying weights, urgency levels, and geographic destinations. To fulfill these order requests, a manager orchestrates delivery cycles following predetermined epochs, distributing packages from the depot via drones. At each decision epoch, the manager must select a subset of orders for immediate fulfillment through a coordinated wave of drone delivery. Upon the completion of each delivery cycle, fulfilled orders are cleared from the system, whereas unfulfilled ones join with newly incoming requests for consideration in the subsequent cycle. The primary challenge within this logistical framework is the adept management of dynamic and stochastic delivery orders. Rather than simply arranging delivery plans based on existing requests, more efficient strategies should account for potential order arrivals in the near future. Capitalizing on opportunities to batch orders arriving at different times and delivering them collectively can yield significant savings in vehicle miles traveled, particularly when their destinations are close in distance. But to facilitate order consolidation, certain requests may need to be strategically postponed to future delivery cycles, which are as yet unrealized and subject to uncertainty. However, by leveraging historical service records, the manager possesses insights into potential demand profiles. Consequently, the manager must weigh the trade-off between existing orders in the pipeline and potential future arrivals to determine whether to prioritize or postpone individual orders with varying penalties and locations, constrained by the delivery capabilities of drones in both the current and future cycles. This decision-making process necessitates proactive decisions at each epoch and the implementation of corrective actions through recourse rules in subsequent cycles.

Our paper introduces a stochastic optimization model with a two-stage structure designed to address the aforementioned cyclic order fulfillment decisions. The first stage models the decisions and system operations within the current delivery cycle, while the second stage anticipates operations in future cycles, dependent on decisions from the first stage, and feeds potential costs back to the first-stage objective. Through repeated implementation of this model from one cycle to the next, the order fulfillment problem is solved dynamically to minimize total operational costs accrued over time. The main contributions of the paper are summarized as follows:

\begin{itemize}

\item A two-stage stochastic model is formulated to tackle the dynamic order fulfillment problem in on-demand last-mile delivery with drones, taking into account various sources of demand uncertainty, including the variations in arrival patterns, delivery destinations, and package weights of individual orders.

\item An integer L-shaped method is customized to enhance its computational efficiency on the proposed two-stage model. By streamlining the branch-and-cut process to focus on a slim set of variables connecting the first- and second-stage problems, we notably reduce the size of the branching tree. Our experiments demonstrate that this `slim' integer L-shaped method significantly outperforms the brute-force solution of the extended form by Gurobi, which fails to reach optimality within a reasonable timeframe.

\item The second-stage subproblem, which involves routing constraints for drones, acts as a bottleneck in the computation time of the entire problem. In order to meet the requirements of online decision-making, the branch-and-cut process is further refined by substituting the exact method for solving the subproblem with a pair of heuristic methods and integrating them alternately. Incorporating these heuristics accelerates the processing speeds by over 20 times compared to the exact L-shaped method, while still maintaining an optimality gap of less than 1\% in our experiments.

\item A comprehensive sensitivity analysis evaluates the potential of deploying the two-stage stochastic model for drone delivery compared to its myopic single-stage counterpart in a dynamic market setting. The significance of accounting for uncertain future demand profiles depends heavily on the sources and levels of the uncertainty. Overall, the stochastic program typically outperforms the single-stage approach by 20\% in operational costs over the long term when uncertainty regarding future demand remains low. The superiority of the stochastic program remains fairly robust against uncertainties in delivery destinations and package weights. It, however, diminishes with increasing uncertainties in order arrivals.

\end{itemize}

The remainder of this paper is organized as follows. \Cref{sec:liter_review} reviews the related work, summarizing previous studies on three pillars relevant to the current study. \Cref{sec:problem_formulation} introduces our problem setting of dispatching drones cyclically for spatially distributed tasks and formulates a two-stage stochastic program for tackling the problem with intrinsic demand uncertainty. \Cref{sec:int_L-shaped_method,sec:computational_analysis} details the integer L-shaped method and our modifications, followed by numerical results comparing the computational performance of three variants of our customized method to the standard integer L-shaped method and the brute-force solution by Gurobi. \Cref{sec:operational_analysis} conducts further numerical analysis on a series of representative instances to systematically evaluate the benefits of incorporating uncertain future demand in the operations of on-demand drone delivery. At the end, \Cref{sec:conclusion} concludes the paper.

\section{Literature Review}\label{sec:liter_review}
This section reviews the literature along three relevant areas and delineates how our work differs from and contributes to the existing discussions. Firstly, our focus lies in the application of drones in last-mile delivery, distinguished by their distinct and heightened operational constraints compared to conventional delivery vehicles, highlighting the need for the development of sophisticated algorithms for dynamic order fulfillment and routing with drones. Subsequently, in the second area, we delve into existing research related to dynamic and stochastic vehicle routing problems, and further motivate our modeling of dynamic order fulfillment as a two-stage stochastic program with rolling horizons. Finally, we provide a succinct review of various modeling and decomposition techniques within the realm of two-stage stochastic integer programming.

\subsection{Last-Mile Delivery using Drones}

The continuous growth of last-mile delivery has been accompanied by various business innovations, including the testing of autonomous delivery means to enhance direct-to-customer services \citep{straits2023global}. Among these solutions, there are three main forms of autonomous delivery vehicles (ADVs) covered: drones or UAVs, sidewalk autonomous delivery robots, and road autonomous delivery robots \citep{figliozzi2020carbon}. The use of these ADVs has been proof-concepted to benefit the last-mile delivery industry in different aspects of operations \citep{alverhed2024autonomous}, ranging from lower costs \citep{hoffmann2018regulatory, boysen2018scheduling, heimfarth2022mixed}, higher efficiency \citep{jennings2019study}, and reduced carbon footprint \citep{figliozzi2020carbon}. 

Compared to ground-based ADVs, delivery drones have lower carrying capacities but achieve significantly higher speeds and greater flexibility by bypassing complex road conditions and traveling along Euclidean paths. This makes them particularly suitable for long-distance and time-sensitive deliveries. Research on drone-enabled delivery operations has a long history and is primarily categorized into two main modeling approaches: drone-only operations and truck-drone collaborative systems. The drone-only models are often treated as variants of the classic vehicle routing problem (VRP) but recognize specific limitations and constraints for drone operations. Particularly, due to limited battery capacity, package delivery with drones needs to take into account battery charging and degradation into lifecycle cost analysis \citep{liu2023routing}. Unlike ground vehicles, the energy consumption of drones depends on both payload and travel distance \citep{cheng2020drone}. Overlooking load-dependent constraints may underestimate travel costs and may potentially lead to infeasible routing for drones \citep{xia2023branch}. Other unique considerations for drone delivery include limited carrying capacity \citep{dorling2016vehicle}, no-fly zones \citep{jeong2019truck}, safety separation assurance \citep{chen2024integrated}, and varying wind conditions \citep{sorbelli2020energy}. Some recent studies have built upon more sophisticated truck-drone tandem models, wherein drones are launched from delivery trucks and coordinate their routing alongside trucks to efficiently deliver packages. A single drone may execute multiple flights and landings, potentially from and to different trucks, while a truck might launch and retrieve several drones at various time and locations \citep{murray2015flying}. For this category of models, determining the service sequence for customers by both vehicle types, pinpointing the locations for drone deployment and collection by trucks, and mapping out the vehicles' routes presents new challenges to the VRPs \citep{wang2019vehicle, tamke2021branch}. This study continues the above thread of discussions by focusing on a drone-only model, having a fleet of drones dispatched from a central depot through delivery cycles for on-demand order fulfillment. Drawing upon the considerations of existing research, each drone has a specified battery capacity that imposes constraints on the travel range per delivery tour. Additionally, the energy consumption rates of drones are determined by their en-route payload and will be optimized in conjunction with other operational costs in the dispatch and touring of drones.

While optimization models for drone delivery systems have received growing attention, problems that incorporate uncertainty in system conditions remain relatively understudied. Existing research has mostly focused on strategic decision-making. For example, to enhance the timely delivery of aid packages to disaster-affected areas, \cite{kim2019stochastic} investigated the optimal placement of drone take-off platforms by incorporating a stochastic flight coverage range influenced by weather conditions. Similarly, \cite{ghelichi2022drone} considered a drone facility location problem that accounts for uncertainty in demand locations. Beyond location decisions, some studies also examined the optimal number of drones while accounting for uncertainties in battery duration \citep{kim2018drone} or demand \citep{hadas2024modeling}. Stochastic scheduling problems have also been explored. For instance, \cite{huang2020reliable} developed a model to maximize the probability of completing a drone journey within a given time threshold on a stochastic, time-dependent public transportation network. Likewise, \cite{deng2024stochastic} proposed a robust optimization model to minimize service delays by incorporating stochastic travel times in a truck-drone collaborative delivery system. In contrast to these studies, which focus on one-time delivery routing and scheduling, our work addresses real-time decision-making in a dynamic environment. The delivery system under consideration operates within a rolling-horizon framework, where the operational decision is to exercise controlled delays in each cycle to enable opportunities for batching orders and reducing costs. Methodologically, this objective aligns more closely with dynamic and stochastic VRPs, where recourse actions are allowed and even encouraged to correct previous delays.

\subsection{Dynamic and Stochastic Vehicle Routing Problem}

Dynamic and stochastic VRP arises in on-demand services such as ride-sharing and meal or medication delivery, where extreme timeliness is essential for navigating a transient and stochastic market environment. A widely used approach to address method of tackling such dynamic decision-making problems is through the application of Markov decision processes (MDPs), which have also been employed to model and address vehicle dispatch and routing for order fulfillment problems \citep{yoon2021dynamic, jenkins2021approximate}. Dynamic programming as a common technique for solving MDPs, however, faces significant challenges in large-scale stochastic environments owing to the explosive increase in the dimensions of states, decisions, and scenario realizations, especially when involving vehicle routing. Consequently, heuristic methods or approximations are often needed for adopting MDPs in such use cases. \cite{secomandi2009reoptimization} developed a finite-horizon MDP formulation for the VRP with stochastic demand and solved it through a partial re-optimization heuristic by restricting all possible states to a small subset. Similarly, \cite{klapp2018one} formulated the dynamic dispatch waves problem as an MDP, where orders are batched through delivery waves. They solved this MDP by approximating the cost-to-go function and converting the MDP into an approximate linear program. Contrary to dynamic programming, which explores the entire state and action space of MDPs, multi-period optimization with rolling horizons takes further approximations by seeking local optima for individual or adjacent decision epochs without backtracking the entire future trajectories. The rolling-horizon framework typically couples with online algorithms, making decisions based on the current state and a limited future horizon \citep{alden1992rolling}. This method has been widely applied to dynamic problems such as energy planning \citep{silvente2015rolling}, resource allocation \citep{bertsimas2017comparison}, and vehicle routing problems \citep{figliozzi2007pricing, ma2023dynamic}. In this paper, we opt for the latter approach to address an online dynamic order fulfillment problem using rolling horizons instead of solving offline a priori policies. To better adapt to the dynamically revealed stochastic demand during execution, we model the optimization for each decision epoch as a stochastic order dispatch and vehicle routing problem.

The formulations for stochastic order fulfillment fall into two categories: chance-constrained programs and stochastic programs with recourse. For the former, delivery failures are permitted, with the probability of failure constrained to remain below a certain threshold \citep{raff1983routing}. The more prevalent formulation is the latter, aiming to find a solution that optimizes objectives considering both certain and uncertain aspects of system realizations. In the context of order fulfillment, uncertainty could stem from the demand side \citep{bertsimas1992vehicle, laporte2002integer} or the supply side\citep{laporte1992vehicle, tacs2013vehicle}. Specifically, \cite{bertsimas1992vehicle} and \cite{laporte2002integer} considered stochastic demand with uncertain package weights. They assumed the weight of each package was unknown before arriving at the customers and vehicles might be unable to load the full customer demand, resulting in return trips to the depot. These return trips are followed by the resumption of the original routes and incur additional penalties due to corrective trips. Our problem entails different forms of uncertainty and managerial goals. All pending orders are fully realized in our case, eliminating the need for routing corrections. Instead, our decision-making revolves around whether to accept pending orders in the current or future cycles.

\subsection{Decomposition Methods for Two-Stage Stochastic Integer Programming}

To solve the stochastic programs with recourse, the most common approach is to convert it to an equivalent deterministic problem by constructing a sample average formulation. The basic idea is to approximate the recourse function through the average estimates from a finite set of scenarios sampled based on known distributions of uncertain events. When the number of scenarios is sufficiently large, the obtained optimal solution will converge to the true optimal value almost surely \citep{birge2011introduction}. However, accurately solving sample average models with large simulated scenarios can be time-consuming. As a result, the implementation of the integer L-shaped method, an exact decomposition algorithm is effective for a broad spectrum of stochastic integer programs that involve recourse. This branch-and-cut method iteratively computes a first-stage solution by introducing additional optimality cuts generated using the immediate expected cost of recourse in each iteration \citep{laporte1993integer}. When compared to directly solving the sample average formulation with a large number of scenarios, the integer L-shaped method has demonstrated remarkable efficiency on many applications, including vehicle routing problem \citep{laporte2002integer}, location and network restoration problem \citep{sanci2021integer} and generalized assignment problem \citep{albareda2006exact}.

As an enhancement to the original integer L-shaped method by \cite{laporte1993integer}, \cite{angulo2016improving} further improved the efficiency through two key strategies of alternating between linear and mixed integer subproblems and introducing a new type of integer optimality cut. Building upon this framework of alternating cuts, \cite{sanci2021integer} took a modification by employing a special branching scheme based on the fractional condition of the variables of interest, as opposed to the traditional binary branching tree. Despite these advancements, the modified integer L-shaped method, aimed at exact solutions, still faced computational challenges that made it unsuitable for real-time decision-making. One alternative solution is the cut-and-project framework developed by \cite{bodur2017strengthened} to strengthen the L-shaped cut, but it requires continuous recourse (second-stage) variables, which are inapplicable in our case. \cite{mousavi2022stochastic} identified a totally unimodular structure in their second-stage problem and applied the cut-and-project framework on a stochastic last-mile delivery problem with uncertain availability of crowd-shippers. Nonetheless, our problem presents a more intricate second-stage structure that spans several delivery cycles, making it difficult to adapt to this framework. The seminal work closely related to ours is by \cite{larsen2023fast}, who introduced a groundbreaking approach by substituting the exact second-stage solutions with fast, yet accurate, predictions from supervised machine learning. Similarly, \cite{chan2022machine} incorporated a machine learning prediction model into a reduced sampling strategy to approximate the objective function of the unsampled second-stage set, addressing bilevel optimization problems involving a large number of independent followers. Instead of pre-training an offline prediction model, we adopt the heuristics provided by Google OR-tools for solving second-stage subproblems. Additionally, we integrate them within a two-tier structure for generating optimality cuts, ultimately enhancing computational performance to handle on-demand operations.

\section{Problem Formulation}\label{sec:problem_formulation}

\begin{figure}[b!]
	\centering
	\includegraphics[width=1.0\textwidth]{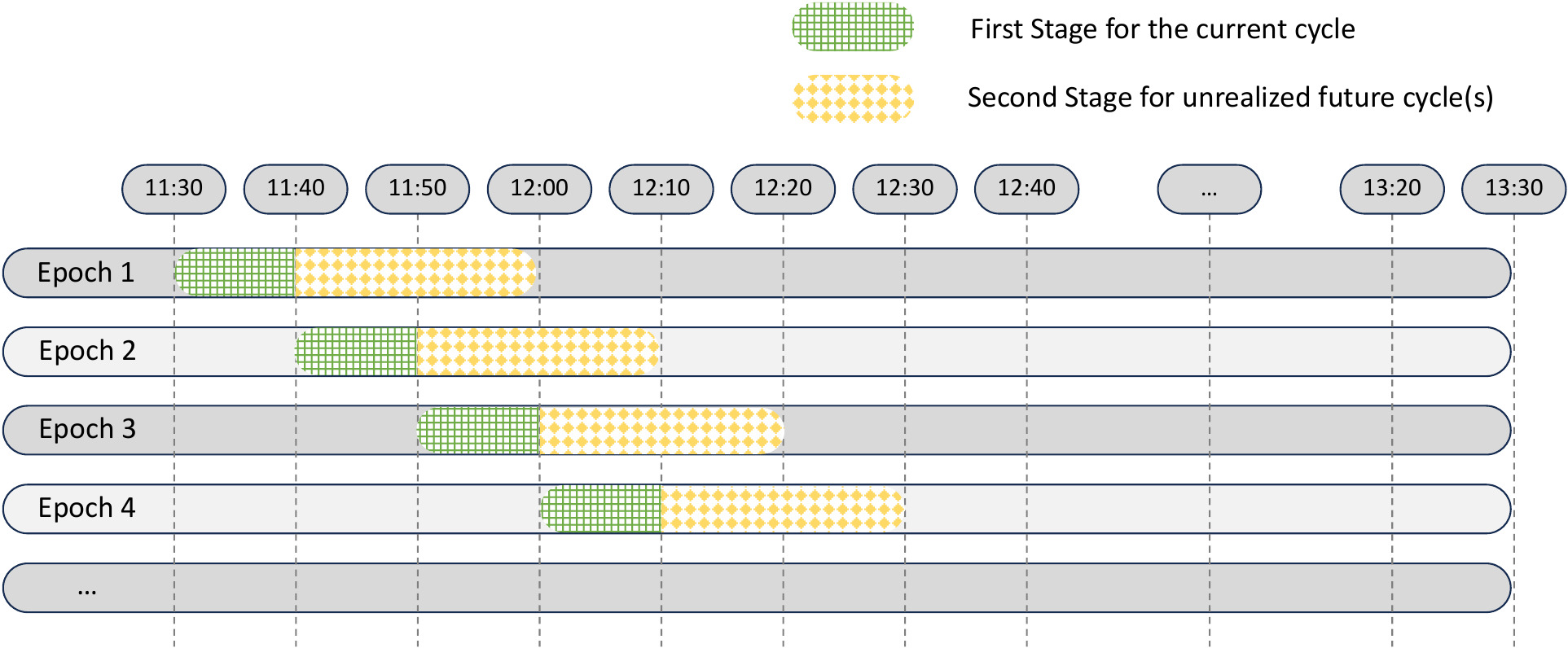}
	\caption{The rolling-horizon implementation of the two-stage stochastic program for dynamic order fulfillment at a service depot. The first stage takes the operational decisions for the current delivery cycle at each epoch, with the second stage accounting for their interactions and consequences across several subsequent cycles within a defined future horizon.}
	\label{rolling_graph}
\end{figure}

This section formulates a two-stage stochastic model for cyclic order fulfillment using drones dispatched from a service depot, illustrated in Figure \ref{rolling_graph}. At each decision epoch, the depot selects a subset of pending requests for immediate fulfillment in the current delivery cycle using a fleet of drones. We assume the depot has access to prospective demand profiles, which encompass distributions of the number of package arrivals, along with their weights, urgency levels, penalties for late delivery, and destinations. Consequently, when optimizing decisions in the current cycle, this future information should be incorporated. Taking the decision Epoch 3 at 12:00 in Figure \ref{rolling_graph} as an example, our proposed stochastic model makes first-stage decisions that involve selecting pending orders actualized at 12:00 and dispatching drones to fulfill them in the current delivery cycle. Any order not accommodated within the current cycle will be deferred to the subsequent cycles at 12:10 and 12:20, incurring heightened penalties. In the meanwhile, the second stage pictures prospective demand arrivals during the period from 12:00 to 12:20 and models the interactions between current fulfillment decisions and prospective system realizations. The ultimate goal of the two-stage model is to facilitate optimal decision-making for the current cycle, aiming to minimize the expected operational cost in the long run.

The model distinguishes between the current and future cycles using the superscript $t$. All states and variables related to the current delivery cycle are associated with $t=0$, while subsequent future cycles within the decision horizon are denoted as $t\in \netT = \{ 1,2,\dots,T \}$, where the index $T$ denotes the cardinality of $\netT$. The delivery cycles are predetermined, and the length of each cycle is assumed sufficient for a drone to travel back and forth from the furthest possible delivery point. Each drone has limited load capacity $Q$ and a battery capacity $E$, with which the dispatching and routing decisions of drones must comply. By the end of each delivery cycle, all dispatched drones will return to the depot to be replenished with fully charged batteries for the subsequent cycle of deliveries. Based on this setup, the remaining sections elaborate on the modeling components separately for the first and second stages, and then integrate both stages into a unified model. A notational glossary is enclosed in Table \ref{notatoin1} at the end of the manuscript for reference.

\subsection{First-Stage Decisions and Constraints}
Consider an arbitrary decision epoch with an accumulation of pending requests, denoted as $M$, which includes both unfulfilled and newly incoming orders from the previous delivery cycle. For each of these existing requests $m \in M$, its exact time of generation, package weight $q_m$, drop-off destination, and time-moving penalties are known. The penalty incurred at any cycle $t(\in \{0\}\cup\netT)$ is predetermined as $c_m^t$. We use the binary variable $y_m^0$ to indicate whether request $m$ will be selected for delivery in the current cycle ($t=0$). 

Additionally, a set of dispatchable drones is available for service, denoted as $K$, and the binary variable $r^{0,k}$ indicates whether drone $k (\in K)$ will be dispatched in cycle $0$. For each drone dispatched in a delivery cycle, it tours through a series of delivery points for package drop-offs. The tours of drones can feasibly traverse from the depot to the drop-off location of an order $m\in M$, from one drop-off location to another, and return back to the depot after dropping off the last package. We assume that within each delivery cycle, each drone executes only one tour that starts and ends at the depot. In other words, drones do not return to the depot and embark on additional drop-off tours within a cycle. To facilitate the modeling of drone routing, we use two aliases to represent the service depot: $O$ and $D$, which respectively denote the starting and returning points for drones heading out for delivery. Thus, mathematically, the set of feasible arcs $A$ for drones can be defined as follows:
\begin{align*}
    A = \{ (i,j) | i \in M \cup \{O\}, j \in M \cup \{D\}, i \neq j \} \backslash \{(O,D)\}.
\end{align*}
To streamline the notation, we abuse the index of order request $m(\in M)$ in an arc tuple to denote the drop-off location of order $m$ serving as the starting or ending node of the arc, respectively. With that, a set of binary decision variables $x_{ij}^{0,k}$ is further introduced to indicate whether drone $k$ traverses through arc $(i,j)$ in cycle $0$, for all $k\in K$ and $(i,j)\in A$.

The selection and dispatch of orders to specific drones as well as the delivery sequence of each drone for assigned packages in the current cycle are subject to the following constraints:
\begin{align}
    & r^{0,k} = \sum_{m \in M} x_{Om}^{0,k} & \forall k \in K, \label{1st-con1} \\
    & r^{0,k} \leq r^{0,k-1} & \forall k \in K \backslash \{1\}, \label{1st-con2} \\
    & y_m^0 = \sum_{k \in K} \sum_{j \in M \cup \{D\}} x_{mj}^{0,k} & \forall m \in M, \label{1st-con3} \\
    & \sum_{i \in M \cup \{O\}} x_{im}^{0,k} - \sum_{j \in M \cup \{D\}} x_{mj}^{0,k} = 0 & \forall m \in M,\ k \in K. \label{1st-con4}
\end{align}
Constraint \eqref{1st-con1} indicates that if drone $k$ is dispatched in the current cycle, it will travel out of the depot to one of the delivery points. Constraint \eqref{1st-con2} is designed for symmetry breaking, ensuring that drone $k$ can be dispatched only when drone $k-1$ is used. Constraint \eqref{1st-con3} guarantees the delivery of order $m$ by one of the drones if the order is selected for fulfillment in the current cycle. Lastly, Constraint \eqref{1st-con4} maintains individual drop-off locations as intermediate visiting points of the drones' delivery tours. 

Furthermore, two continuous variables $u_i^{0,k}$ and $v_i^{0,k}$ are introduced to represent the remaining load of packages (measured in kilograms) and battery power (measured in kilowatt-hours) of drone $k$ immediately after visiting node $i$. Specifically, for any drone $k$ out for delivery, its starting state of charge, $v_O^{0,k}$, at the depot is assumed to be equal to the battery capacity $E$, while the ending load, $u_D^{0,k}$, when returning to the depot is fixed at $0$. The en-route energy consumption rate of drones may vary depending on the load of packages \citep{xia2023branch}. Let function $P(u)$ denote the per-unit distance energy consumption rate of drones with a carrying load of $u$. Then, the energy consumption of one drone traveling from node $i$ to $j$ with load $u_i$ equals to $P( u_i ) \cdot d_{ij}$, where $d_{ij}$ represents the distance traveled between the nodes. Therefore, the reduction of package loads and battery state of charge of each drone through the delivery tours can be captured by the following constraints:
\begin{align}
    & u_j^{0,k} \leq u_i^{0,k} - q_j + Q \cdot (1-x_{ij}^{0,k}) & \forall (i,j) \in A,\ k \in K,\label{1st-con5} \\
    & u_O^{0,k} \leq W,\quad u_D^{0,k} = 0 & \forall k \in K, \label{1st-con6} \\
    & v_j^{0,k} \leq v_i^{0,k} - P \left( u_i^{0,k} \right) \cdot d_{ij} + E \cdot (1-x_{ij}^{0,k}) & \forall (i,j) \in A,\ k \in K, \label{1st-con7} \\
    & v_O^{0,k} = E,\quad v_D^{0,k} \geq 0 & \forall k \in K. \label{1st-con8}
\end{align}

\noindent Constraints \eqref{1st-con5} and \eqref{1st-con7} monitor the reduction of package weights and battery power carried by drones during the delivery tours. Constraints \eqref{1st-con6} and \eqref{1st-con8} ensure that carried packages and energy consumption of each delivery tour do not exceed drones' load capacity and battery capacity, respectively. Together, these four constraints serve as subtour elimination measures \citep{toth2002vehicle}, ensuring a single, coherent route of delivery for each drone.

\subsection{Second-Stage Variables and Constraints}
In the second stage, the model predicts stochastic demand patterns several cycles into the near future, represented by a distribution of scenarios denoted with the superscript $\xi$.

Let $N$ be the set of potential order requests that could arrive within the next $T$ cycles. Each request $n(\in N)$ can realize with different arrival times, delivery locations, and time-moving penalties ${c_n^{t,\xi}}$ in different scenarios. In particular, we let $a_n^{\xi}$ denote the first decision epoch for request $n$ immediately following its arrival in scenario $\xi$. Accordingly, the order selection decisions in the second stage involve both the remaining unfulfilled requests from the current delivery cycle and those new arrivals in the future scenarios. Therefore, we use binary variables $y_m^{t,\xi}$ and $y_n^{t,\xi}$ to indicate whether an existing request $m \in M$ and a newly arriving order $n \in N$ would be fulfilled at the particular cycle $t\in\netT$ in the future. Any orders that cannot be fulfilled within the horizon $\mathcal{T}$ will be automatically considered ``handled'' at an ad-hoc cycle $T+1$. The following three constraints specify the order fulfillment decisions in the second stage as well as their connection to first-stage decisions:
\begin{align}
    & y_m^0 + \sum_{t\in \netT \cup \{T+1\}} y_m^{t,\xi} = 1 & \forall m \in M, \label{2nd-con9} \\
    & \sum_{t\in \netT \cup \{T+1\}} y_n^{t,\xi} = 1 & \forall n \in N, \label{2nd-con10} \\
    & \mathcal{M}\cdot (1-y_n^{t,\xi}) \geq a_n^{\xi} - t & \forall n\in N,\ t\in \netT \cup \{T+1\} \label{2nd-con11},
\end{align}
Specifically, Constraint \eqref{2nd-con9} ensures that existing orders will be fulfilled either in the current cycle or one of the future cycles. Constraint \eqref{2nd-con10} specifies that potential orders will be handled in a future cycle, while the Big-M Constraint \eqref{2nd-con11} restricts that those orders can only be serviced after their generation.

It is worth noting that the number of order arrivals in future scenarios can also be stochastic in nature. As the two-stage formulation does not recognize a variable set of $N$ (i.e., $\left | N \right |$ is fixed for different scenarios in $\xi$), we employ a workaround by expanding the delivery range to bypass this restriction and enable randomness and variations in the number of requests. For instance, if the delivery points of orders are distributed regularly within a range of (0, 10] miles from the depot, we could simulate a fixed number of requests over the range of (-5, 10] miles by including an extended negative range. Then, those cases inside the negative range will be treated as invalid requests and served immediately by a ``virtual drone'' with an unlimited service capacity in each cycle. We introduce an additional variable $z_n^{\xi}$ to indicate the validity of order $n$ and whether it is taken care of by the virtual drone once it is generated. Constraints \eqref{2nd-con12}-\eqref{2nd-con15} dictate the second stage's dispatching of selected orders to drones for fulfillment:
\begin{align}
    & y_m^{t,\xi} = \sum_{k \in K} \sum_{j \in M \cup N \cup \{D\}} x_{mj}^{t,k,\xi} & \forall m \in M,\ t\in \netT, \label{2nd-con12} \\
    & y_n^{t,\xi} = \sum_{k \in K} \sum_{j \in M \cup N \cup \{D\}} x_{nj}^{t,k,\xi} + z_n^{\xi} \cdot \mathds{1}(a_n^\xi=t) & \forall n \in N,\ t\in \netT, \label{2nd-con13} \\ 
    & d_{On}^{\xi} \geq -\mathcal{M} \cdot z_n^{\xi} & \forall n\in N, \label{2nd-con14} \\ 
    & d_{On}^{\xi} \leq \mathcal{M} \cdot \left(1 - z_n^{\xi}\right) & \forall n\in N. \label{2nd-con15}
\end{align}
For an order selected to be fulfilled in cycle $t$, Constraints \eqref{2nd-con12} and \eqref{2nd-con13} ensure it is processed by one drone within the cycle. Invalid orders, as indicated by $z_n^{\xi}$, will be handled by a virtual drone. The latter two Big-M Constraints \eqref{2nd-con14} and \eqref{2nd-con15} recover the validity of second-stage orders based on their generated delivery distances from the depot, i.e., $d_{On}^{\xi}$ for $n\in N$.

The set of feasible arcs for drones in the second stage is denoted as $\tilde{A}$, which can be materialized as follows,
\begin{align*}
    \tilde{A} = \{ (i,j) | i \in M \cup N \cup \{O\}, j \in M \cup N \cup \{D\}, i \neq j \} \backslash \{(O,D)\}.
\end{align*}
Similar to the first stage, the delivery tours of drones at each cycle $t \in \netT$ in the second stage should also satisfy the feasibility constraints as follows,
\begin{align}
    & r^{t,k,\xi} = \sum_{j \in M \cup N} x_{Oj}^{t,k,\xi} & \forall k \in K ,\ t\in \netT, \label{2nd-con16} \\
    & r^{t,k,\xi} \leq r^{t,k-1,\xi} & \forall k \in K \backslash \{1\},\ t\in \netT, \label{2nd-con17} \\
    & \sum_{i \in M \cup N \cup \{O\}} x_{ih}^{t,k,\xi} - \sum_{j \in M \cup N \cup \{D\}} x_{hj}^{t,k,\xi} = 0 & \forall h \in M \cup N,\ k \in K,\ t\in \netT, \label{2nd-con18} \\
    & u_j^{t,k,\xi} \leq u_i^{t,k,\xi} - q_j^\xi + Q \cdot (1-x_{ij}^{t,k,\xi}) & \forall (i,j) \in \tilde{A},\ k \in K,\ t\in \netT, \label{2nd-con19} \\
    & u_O^{t,k,\xi} \leq Q,\ u_D^{t,k,\xi} = 0 & \forall k \in K,\ t\in \netT, \label{2nd-con20} \\
    & v_j^{t,k,\xi} \leq v_i^{t,k,\xi} - P \left( u_i^{t,k,\xi} \right) \cdot d_{ij}^\xi + E \cdot (1-x_{ij}^{t,k,\xi}) & \forall (i,j) \in \tilde{A},\ k \in K,\ t\in \netT, \label{2nd-con21} \\
    & v_O^{t,k,\xi} = E,\ v_D^{t,k,\xi} \geq 0 & \forall k \in K,\ t\in \netT. \label{2nd-con22}
\end{align}
Again, Constraints \eqref{2nd-con16} and \eqref{2nd-con17} are for symmetric-breaking dispatching of drones. Constraint \eqref{2nd-con18} maintains each delivery point as an intermediate visiting node of delivery tours. Constraints \eqref{2nd-con19}-\eqref{2nd-con22} track the en-route changes in drones' package loads and battery power and eliminate subtours for delivery.

\subsection{Two-Stage Stochastic Program}
Let $\bm{\netX^0}$ represent the tuple summarizing all the first-stage decision and state variables, denoted as $\bm{\netX^0} = \{\bm{x^0, y^0, r^0, u^0, v^0}\}$. Below, we massage the aforementioned two stages of constraints into a unified program for optimizing the order fulfillment decisions at the current epoch, in consideration of stochastic future demand:

(\textbf{P1}):
\begin{subequations}\label{opt_model_1}
\begin{align}
    & \min_{\bm{x^0, y^0, r^0} \in \mathds{B}; \bm{u^0, v^0} \in \mathds{R^+}}\ \alpha \cdot \sum_{m \in M} c_m^0 \cdot y_m^0 + \beta \cdot \sum_{(i,j) \in A} \left( d_{ij}\cdot \sum_{k \in K} x_{ij}^{0,k}\right) + \gamma \cdot \sum_{k \in K} \left(v_O^{0,k} - v_D^{0,k}\right) \nonumber \\
    & \qquad \qquad \qquad \ \  + \kappa \cdot \sum_{k \in K} r^{0,k} + \mathbb{E}_{\xi} \left[ \mathbb{Q}(\bm{\netX^0}, \xi) \right] \\
    & \text{s.t.} \quad \quad (\ref{1st-con1}) - (\ref{1st-con8}) \nonumber \\
    & \text{where} \nonumber \\
    & \mathbb{Q}(\bm{\netX^0}, \xi) = \min_{\bm{x^{t,\xi}, y^{t,\xi}, z^{t,\xi}, r^{t,\xi} \in \mathds{B}}; \bm{u^{t,\xi}, v^{t,\xi}} \in \mathds{R^+}} \ \alpha \cdot \sum_{t\in \netT \cup \{T+1\}} \left( \sum_{m \in M} c_m^t \cdot y_m^{t,\xi} + \sum_{n \in N} c_n^{t,\xi} \cdot  y_n^{t,\xi} \right) \nonumber \\
    &+ \sum_{t\in \netT} \left( \beta \cdot \sum_{(i,j) \in \tilde{A}} \left( d_{ij}^{\xi}\cdot \sum_{k \in K} x_{ij}^{t,k,\xi}\right) + \gamma \cdot \sum_{k \in K} \left(v_O^{t,k,\xi}-v_D^{t,k,\xi}\right) + \kappa \cdot \sum_{k \in K} r^{t,k,\xi} \right)  \\
    & \text{s.t.} \quad \quad (\ref{2nd-con9}) - (\ref{2nd-con22}). \nonumber
\end{align}
\end{subequations}
In this optimization problem, the objective function (\ref{opt_model_1}a) seeks to minimize the total operational costs incurred in the current cycle and the expected costs in future cycles, by taking the expectation of the second-stage value function in (\ref{opt_model_1}b) over $\xi$. Both the first- and second-stage objectives comprise the following components in order: the first term corresponds to the penalty for delayed order delivery; the second term pertains to variable costs of drones that are proportional to touring distances; the third term addresses the cost for drones' energy consumption; and the final term accounts for the fixed cost associated with each time of use for drones. The coefficients for these cost terms, denoted as $\alpha$, $\beta$, $\gamma$, and $\kappa$, are constants with values specified exogenously.

\subsection{Sample Average Approximation}

Solving the above optimization problem (\textbf{P1}) requires the computation of the expected values outputted by the second-stage value function $\mathbb{Q}(\bm{\netX^0}, \xi)$. In nature, this characterizes an integral in extremely high dimensions that is challenging, if not impossible, to compute. The most common solution is to solve an equivalent deterministic sample-average problem by sampling a discrete set of realizations $\Omega$ instead of continuously distributed random events $\xi$ \citep{birge2011introduction}. Then, instead of computing the expectation of $\mathbb{Q}(\bm{\netX^0}, \xi)$ in problem (\textbf{P1}), the sample average approximation (SAA) is used as a surrogate based on a finite set of scenarios $\Omega$.

Following this path, this subsection reformulates (\textbf{P1}) to the extensive form (\textbf{EF}) of SAA as an integer program:

(\textbf{EF}):
\begin{subequations}\label{opt_model_2}
    \begin{align}
        & \min_{\bm{x, y, r, z} \in \mathds{B}; \bm{u, v} \in \mathds{R^+}}\ \alpha \cdot \sum_{m \in M} c_m^0 \cdot y_m^0 + \beta \cdot \sum_{(i,j) \in A} \left( d_{ij}\cdot \sum_{k \in K} x_{ij}^{0,k}\right) + \gamma \cdot \sum_{k \in K} \left(v_O^{0,k} - v_D^{0,k}\right) + \kappa \cdot \sum_{k \in K} r^{0,k} \nonumber \\
        & + \frac{1}{\left | \Omega \right |} \sum_{\omega \in \Omega} \left\{ \alpha \cdot \sum_{t\in \netT \cup \{T+1\}} \left( \sum_{m \in M} c_m^t \cdot y_m^{t,\omega} + \sum_{n \in N} c_n^{t,\omega} \cdot  y_n^{t,\omega} \right) \nonumber \right.\\
        &\quad\quad + \sum_{t\in \netT} \left. \left( \beta \cdot \sum_{(i,j) \in \tilde{A}} \left(d_{ij}^{\omega}\cdot \sum_{k \in K} x_{ij}^{t,k,\omega}\right) + \gamma \cdot \sum_{k \in K} \left(v_O^{t,k,\omega}-v_D^{t,k,\omega}\right) + \kappa \cdot \sum_{k \in K} r^{t,k,\omega} \right) \right\} 
    \end{align}
    \begin{align}
        & \text{s.t.} \quad \quad (\ref{1st-con1}) - (\ref{1st-con8}) \nonumber \\
        & y_m^0 + \sum_{t\in \netT \cup \{T+1\}} y_m^{t,\omega} = 1 & \forall m \in M, \omega \in \Omega \\
        & \sum_{t\in \netT \cup \{T+1\}} y_n^{t,\omega} = 1 & \forall n \in N, \omega \in \Omega \\
        & \mathcal{M} \cdot (1-y_n^{t,\omega}) \geq a_n^{\omega} - t & \forall n\in N, t\in \netT \cup \{T+1\}, \omega \in \Omega \\
        & y_m^{t,\omega} = \sum_{k \in K} \sum_{j \in M \cup N \cup \{D\}} x_{mj}^{t,k,\omega} & \forall m \in M, t\in \netT, \omega \in \Omega \\
        & y_n^{t,\omega} = \sum_{k \in K} \sum_{j \in M \cup N \cup \{D\}} x_{nj}^{t,k,\omega} + z_n^{\omega} \cdot \mathds{1}(a_n^\omega=t) & \forall n \in N, t\in \netT, \omega \in \Omega \\ 
        & d_{On}^{\omega} \geq -\mathcal{M} \cdot z_n^{\omega} & \forall n\in N, \omega \in \Omega \\ 
        & d_{On}^{\omega} \leq \mathcal{M} \cdot \left(1 - z_n^{\omega}\right) & \forall n\in N, \omega \in \Omega \\
        & r^{t,k,\omega} = \sum_{j \in M \cup N} x_{Oj}^{t,k,\omega} & \forall k \in K , t\in \netT, \omega \in \Omega \\
        & r^{t,k,\omega} \leq r^{t,k-1,\omega} & \forall k \in K \backslash \{1\}, t\in \netT, \omega \in \Omega \\
        & \sum_{i \in M \cup N \cup \{O\}} x_{ih}^{t,k,\omega} - \sum_{j \in M \cup N \cup \{D\}} x_{hj}^{t,k,\omega} = 0 & \forall h \in M \cup N, k \in K, t\in \netT, \omega \in \Omega \\
        & u_j^{t,k,\omega} \leq u_i^{t,k,\omega} - q_j^\omega + Q \cdot (1-x_{ij}^{t,k,\omega}) & \forall (i,j) \in \tilde{A}, k \in K, t\in \netT, \omega \in \Omega\\
        & u_O^{t,k,\omega} \leq Q,\ u_D^{t,k,\omega} = 0 & \forall k \in K,\ t\in \netT, \omega \in \Omega\\
        & v_j^{t,k,\omega} \leq v_i^{t,k,\omega} - P \left( u_i^{t,k,\omega} \right) \cdot d_{ij} + E \cdot (1-x_{ij}^{t,k,\omega}) & \forall (i,j) \in \tilde{A}, k \in K, t\in \netT, \omega \in \Omega\\
        & v_O^{t,k,\omega} = E,\ v_D^{t,k,\omega} \geq 0 & \forall k \in K,\ t\in \netT, \omega \in \Omega.
    \end{align}
\end{subequations}
The scenarios in (\textbf{EF}) are generated based on Monte Carlo simulation, which will be detailed later in Section \ref{data_generation}. While (\textbf{EF}) is directly solvable, it is expected to be time-consuming due to the drone routing components involved in the second stage. Thus, an integer L-shaped method is further customized to expedite the solution for our problem.

\section{Integer L-Shaped Method}\label{sec:int_L-shaped_method}

This section delineates the development of a customized integer L-shaped method, specifically designed to address dynamic order fulfillment decisions with the SAA formulation. Our exposition is structured to provide readers with a step-by-step understanding of the methodology. Firstly, we briefly overview the L-shaped method in \Cref{sec:general_L}, introducing the concepts of feasibility and optimality cuts, and then bring up the `slim' integer L-shaped cuts in \Cref{sec_cus_L}, refined to our specific problem structure. However, even with these streamlined cuts, exact methods prove overly burdensome and impractical for real-time decision-making. In response, we leverage the capabilities of Google OR-Tools to generate L-shaped cuts using heuristic approaches. This leads to the establishment of a two-tiered framework, introduced in \Cref{sec_huristic cut}, strategically designed to balance the trade-offs between the quality and computational efficiency of cuts generated. Finally, we present the complete customization of our L-shaped algorithm in \Cref{sec:cus_L_algo}.

\subsection{General Framework of the L-Shaped Method} \label{sec:general_L}

Before heading to our customization, let us first briefly introduce the standard L-shaped method based on a general two-stage stochastic program:

\begin{align}\label{sp_original form}
    \min_{\bm{x}} \{ \bm{c}^T \bm{x}+\E_\xi \left[ \mathbb{Q} (\bm{x},\xi) \right] | A\bm{x}=\bm{b}, \bm{x} \geq 0 \},
\end{align}
where $\bm{x}$ denotes the first-stage decision variables; $\bm{c}^T \bm{x}$ calculates the associated costs with certainty; and $A\bm{x}=\bm{b}$ represents a set of linear constraints imposed on $\bm{x}$. Given the separable structure between the first and second stages, the L-shaped method aims to build an outer linearization of the recourse function, depicting the behavior of $\mathbb{E}_\xi \left[\mathbb{Q} (\bm{x},\xi)\right]$ in response to the ex-ante decision $\bm{x}$.

A new variable $\theta$ is introduced to generate a series of lower bound cutting planes to progressively approximate the recourse function $\mathbb{E}_\xi \left[\mathbb{Q} (\bm{x},\xi)\right]$, denoted as $\mathbb{L}(\bm{x})$, which transform the general form \eqref{sp_original form} into (\textbf{P2}) as follows \citep{birge2011introduction}:

(\textbf{P2}):
\begin{subequations}\label{master_form}
	\begin{align}
		& \min z = \bm{c}^T \bm{x} + \theta \\
		\text{s.t.} & \nonumber\\
		& Ax=b \\
		& D_l \bm{x} \geq d_l & l=1,2,\cdots,g \\
		& E_l \bm{x} + \theta \geq e_l & l=1,2,\cdots,h \\
		& \bm{x} \geq 0,
	\end{align}
\end{subequations}
The above (\textbf{P2}) is called the master problem. It solves to find a proposal first-stage solution $\bm{x}^\eta$ in $\eta^{\text{th}}$ iteration, which is then sent to the second stage for calculating two types of cuts added back to the master problem. These two types of cuts are added sequentially: (i) feasibility cuts (\ref{master_form}c) ensure the feasibility of the second-stage problem with the first-stage solution produced, and (ii) optimality cuts (\ref{master_form}d) provide the linear approximations of $\mathbb{L}(\bm{x})$. In expression (\textbf{P2}), there are $g$ feasibility cuts and $h$ optimality cuts, for which the parameters $\{D_l, d_l, E_l, e_l\}$ involved are computed from the optimal solutions associated with the second-stage problem. For our problem (\textbf{EF}), the second-stage problem is always feasible with any solution $\bm{x}$ that is feasible to the first stage (i.e., orders may be fulfilled in $T$ cycles, or they may incur a significant penalty in the ad-hoc cycle $T+1$). Thus, there is no need to consider feasibility cuts in our problem; only optimality cuts are relevant. Through the successive addition of optimality cuts, the variable $\theta$ increases to better approximate $\mathbb{L}(\bm{x})$ until an optimal solution is found,  indicating that $\mathbb{L}(\bm{x}^*) = \theta^*$. The generation of specific integer optimality cuts will be covered with details in the next subsection.

\subsection{`Slim' Integer L-Shaped Cut}\label{sec_cus_L}
Adapting the integer L-shaped method developed by \cite{laporte1993integer}, we present a branch-and-cut algorithm to decompose the problem for an effective solution. The integer L-shaped method couples the branching scheme commonly employed in integer programming with the concept of optimality cuts from the L-shaped method. Each iteration $\eta$ begins with the solution of the master problem. If the solution of the first-stage binary decision variables $\bm{\netY^0} = \{\bm{x^0, y^0, r^0}\}$ contains fractional numbers, we create two new branches on one of the fractional variables following the usual branch and bound procedure. Otherwise, if the solution satisfies integrality, we need to check whether the condition $\theta^\eta \geq \mathbb{L}(\bm{\netX^{0,\eta}})$ holds. If it does not hold, the following optimality cut, developed by \cite{laporte1993integer}, would be added to the master problem:
\begin{align}\label{opt_cut}
     \theta \geq (\mathbb{L}(\bm{\netX^{0,\eta}})-L) \left(\sum_{i \in S^\eta} \netY_i^0 - \sum_{i \notin S^\eta} \netY_i^0 \right) - (\mathbb{L}(\bm{\netX^{0,\eta}})-L)(\left| S^\eta \right| -1) + L.
\end{align}
For this cut added, $S^\eta$ is a set of indices defined as $S^\eta := \{i: \netY_i^{0,\eta}==1\}$; $\left| S^\eta \right|$ represents the count of first-stage binary variables with a value of one in this iteration; $\mathbb{L}(\bm{\netX^{0,\eta}})$ denotes the corresponding second-stage expected value; $L$ is the lower bound for $\mathbb{L}(\bm{\netX^0})$, i.e., $L \leq \min_{\bm{\netX^0}} \{ \mathbb{L}(\bm{\netX^0}) \ | \ (\ref{1st-con1}) - (\ref{1st-con8})\}$. 

Different from the original treatment proposed by \cite{laporte1993integer}, one modification we make here is that instead of considering all first-stage binary variables $\bm{\netY^0}$, we slim down to the terms associated with $y_m^0$ in constructing the cut. This refinement allows us to build the branching tree around $y_m^0$, but ignore those for $r^{0,k}$ and $x_{ij}^{0,k}$, i.e.,
\begin{align}\label{opt_cut_1}
     \theta \geq (\mathbb{L}(\bm{y^{0,\eta}})-L) \left(\sum_{m \in S^\eta} y_m^0 - \sum_{m \notin S^\eta} y_m^0 \right) - (\mathbb{L}(\bm{y^{0,\eta}})-L)(\left| S^\eta \right| -1) + L,
\end{align}
where $y_m^{0,\eta}=1$ for $m \in S^\eta$ and $y_m^{0,\eta}=0$ for $m \notin S^\eta$ represent the the first-stage decisions on whether to fulfill order $m$ in the current cycle. The following proposition formally states Eq. (\ref{opt_cut_1}) as a valid optimality cut. (The proofs of all propositions are enclosed in Appendix \ref{sec:appd_proof} for reference.)
\begin{prop}\label{proposition_1}
    Eq. (\ref{opt_cut_1}) represents a valid optimality cut for the problem (\textbf{EF}).
\end{prop}
\noindent The rationale behind \Cref{proposition_1} is that in our problem, only the pending requests from the first stage affect the outcomes of the recourse. The distinct routing taken by drones to fulfill the same set of orders within a single cycle would not interact with decisions in future cycles, given that the battery is replenished, ensuring identical starting states for drones at the commencement of each cycle. Therefore, unlike the original cut proposed by \cite{laporte1993integer}, the second-stage value function $\mathbb{Q}(\bm{y^{0}}, \xi)$ as well as its expectation function $\mathbb{L}(\bm{y^0})$ in our case depends solely on $\bm{y^0}$ rather than $\bm{r^0}$ and $\bm{x^0}$. One point worthy of attention is that since the branching tree is not constructed using all first-stage binary variables, we need to recheck the feasibility of the problem under the current solution of $\bm{y^0}$, i.e., the existence of feasible solutions with $\bm{r^0}$ and $\bm{x^0}$ being non-negative and integral under $\bm{y^0}$. For \Cref{cus_algorithm} proposed in \Cref{sec:cus_L_algo}, Step \ref{algo_step4} is designed specifically for such purpose.

By incorporating the modified cut, we formulate the master problem for (\textbf{EF}) as follows:

(\textbf{MP}):
\begin{subequations}\label{MP}
	\begin{align}
            & \min_{\bm{x^0, y^0, r^0}\in [0,1]; \bm{u^0, v^0}, \theta \in \mathds{R^+}}\ \alpha \cdot \sum_{m \in M} c_m^0 \cdot y_m^0 + \beta \cdot \sum_{(i,j) \in A} \left(d_{ij}\cdot \sum_{k \in K} x_{ij}^{0,k} \right) + \gamma \cdot \sum_{k \in K} \left(v_O^{0,k} - v_D^{0,k}\right) \nonumber\\
            & \quad\quad\quad\quad\quad\quad\quad\quad\quad\quad\quad\quad\quad\quad\quad + \kappa \cdot \sum_{k \in K} r^{0,k} + \theta
        \end{align}
	\begin{align}
		\text{s.t.} \quad & (\ref{1st-con1}) - (\ref{1st-con8}) \nonumber\\
            & \theta \geq (\mathbb{L}(\bm{y^{0,\eta}})-L) \left(\sum_{m \in S^\eta} y_m^0 - \sum_{m \notin S^\eta} y_m^0 \right) - (\mathbb{L}(\bm{y^{0,\eta}})-L)(\left| S^\eta \right| -1) + L & \forall \eta \in \text{\textEta}.
	\end{align}
\end{subequations}
In this problem, the binary variables $\bm{x^0}$, $\bm{y^0}$, and $\bm{r^0}$ are relaxed to be continuous as per the branch-and-cut framework for the integer L-shaped method \citep{laporte1993integer}. The term $\mathbb{L}(\bm{y^{0,\eta}})$ calculates the sample average of the second-stage value given $\bm{y^{0,\eta}}$, i.e., $\frac{1}{|\Omega|}\sum_{\omega \in \Omega} \mathbb{Q}(\bm{y^{0,\eta}}, \omega)$. Each $\mathbb{Q}(\bm{y^{0,\eta}})$ solves an instance of subproblem (\textbf{SP}) as follows, which can be processed in parallel for $|\Omega|$ instances: 

(\textbf{SP}):
\begin{subequations}\label{SP}
    \allowdisplaybreaks
	\begin{align}
            & \mathbb{Q}(\bm{y^{0,\eta}}, \omega) = \min_{\bm{x^t, y^t, z^t, r^t \in \mathds{B}}; \bm{u^t, v^t} \in \mathds{R^+}} \ \alpha \cdot \sum_{t\in \netT \cup \{T+1\}} \left( \sum_{m \in M} c_m^t \cdot y_m^{t,\omega} + \sum_{n \in N} c_n^{t,\omega} \cdot  y_n^{t,\omega} \right) \nonumber \\
            &+ \sum_{t\in \netT} \left[ \beta \cdot \sum_{(i,j) \in \tilde{A}} \left(d_{ij}^{\omega} \cdot\sum_{k \in K} x_{ij}^{t,k,\omega}\right) + \gamma \cdot \sum_{k \in K} \left(v_O^{t,k,\omega}-v_D^{t,k,\omega}\right) + \kappa \cdot \sum_{k \in K} r^{t,k,\omega} \right]
        \end{align}
        \begin{align}
            \text{s.t.} \quad & (\ref{opt_model_2}c) - (\ref{opt_model_2}o) \nonumber \\
            & y_m^{0,\eta} + \sum_{t\in \netT \cup \{T+1\}} y_m^{t,\omega} = 1 & \forall m \in M
	\end{align}
\end{subequations}

\subsection{Heuristic Solution for Subproblem (SP)}\label{sec_huristic cut}

Each instance of the subproblem (\textbf{SP}) encapsulates a time-expanded order dispatching and vehicle routing problem, whose complexity escalates significantly with the influx of order requests. Despite leveraging the customized L-shaped method, computing the optimal responses of each scenario realization remains arduous. To mitigate the time-intensive challenges associated with solving the subproblem, \cite{larsen2023fast} advocates for substituting the exact resource function values $\mathbb{Q}(\bm{y^{0,\eta}}, \omega)$ with some high-quality proxies that can be obtained efficiently. They trained supervised machine learning models offline to predict $\mathbb{Q}(\bm{y^{0,\eta}}, \omega)$ and then apply those models for quick predictions as surrogates in tasks requiring speedy decision-making. Their extensive experimental results demonstrate that this substitution strategy considerably reduce online solution times while maintaining negligible optimality gaps in most cases. Inspired by their methodology, we propose to solve (\textbf{SP}) through efficient heuristics, using built-in solvers provided by Google OR-Tool's vehicle routing package, and then take the estimated recourse costs $\mathbb{L}^{heu}(\bm{y^{0,\eta}})$ for generating integer L-shaped cuts. In consonance with the specifications of the OR-Tool API, for each realized scenario $\omega \in \Omega$, we transform the subproblem (\textbf{SP}) as per the framework of a capacitated VRP with time-window constraints for effective solutions. The details of how we synchronize the subproblem with the API are discussed in Appendix \ref{sec:appd_subproblem_acc}.

Specifically, Google OR-Tool offers two different heuristic capabilities for solving VRP. The first approach is a greedy algorithm, which, by design, swiftly generates a feasible but potentially suboptimal solution. In contrast, the second approach adopts a more refined approach, meticulously designed to escape local minima and continue the search for improved solutions. Further elaboration on these heuristics, as well as their comparison to the exact method, is also provided in Appendix \ref{sec:appd_subproblem_acc}. Both of the two heuristics can be used to solve the second-stage subproblem and output the recourse for different scenarios, which are then substituted into \eqref{opt_cut_1} to generate the optimality cuts throughout iterations. To clearly distinguish between them, we introduce the notation $\mathbb{L}^{gred}(\bm{y^{0,\eta}})$ and $\mathbb{L}^{aug}(\bm{y^{0,\eta}})$ to represent the average recourse costs yielded through greedy and augmented heuristics, respectively. Consequently, the derived optimality cuts from them are termed the `greedy cut' and `augmented cut'. Regardless of whether the second-stage problem is solved using exact or any heuristic methods, the lower bound $L$ in the optimality cuts \eqref{opt_cut_1} remains the same, specified as $\mathbb{L}(\bm{y^{0}}=\bm{1})$. This is reasoned by \Cref{proposition_2} below:
\begin{prop}\label{proposition_2}
    For the subproblem \eqref{SP}, $\mathbb{L}(\bm{y^{0}}=\bm{1})$ poses a lower bound on the second-stage recourse costs $\mathbb{L}(\bm{y^{0}})$ for all possible $\bm{y^{0}}$.
\end{prop}

Between the two heuristic optimality cuts, the `greedy cut' is characterized by its computational expediency relative to the `augmented cut'. However, this efficiency may compromise the quality of the resultant solutions, which, despite being generally satisfactory, may not rival the higher-quality solutions produced by the `augmented cut'. Then, because $\mathbb{L}^{gred}(\bm{y}^0) \geq \mathbb{L}^{aug}(\bm{y}^0) \approxeq \mathbb{L}(\bm{y}^0)$ for a minimization problem, replacing the cut $\theta \geq \mathbb{L}^{aug}(\bm{y^{0,\eta}})$ with $\theta \geq \mathbb{L}^{gred}(\bm{y^{0,\eta}})$ could result in overestimation of the recourse costs, potentially leading to the exclusion of optimal recourse decisions and reduced quality of solutions to the two-stage stochastic program. To balance the influence of the two heuristics and leverage their relative strengths, a hyperparameter $\nu$ is introduced to assist in determining whether it is necessary to proceed with $\mathbb{L}^{aug}(\bm{y^{0,\eta}})$ in addition to $\mathbb{L}^{gred}(\bm{y^{0,\eta}})$.
\begin{prop}\label{proposition_3}
    Suppose there exist $\nu\in (0,1)$ such that $\nu \mathbb{L}^{gred}(\bm{y^{0}}) < \mathbb{L}(\bm{y^{0}})$ holds in general. Then, if the condition $\theta^{\eta} < \nu \mathbb{L}^{gred}(\bm{y^{0,\eta}})$ is met in a particular iteration $\eta$, the `augmented cut' is tighter than the `greedy cut' in this case, which means adding the `greedy cut' beforehand will not exclude the optimal solution that would obtained by applying `augmented cut'.
\end{prop}

Inspired by the alternating algorithm presented by \cite{angulo2016improving}, which first implements an easily obtained continuous L-shaped mono-cut before transitioning to the exact integer L-shaped cut, we have devised a similar sequential approximation scheme. In our approach, $\mathbb{L}^{gred}(\bm{y^{0,\eta}})$ serves as an initial proxy for $\mathbb{L}(\bm{y^{0,\eta}})$ for each iteration $\eta$. Based on Proposition \ref{proposition_3}, we then verify if the condition $\theta^{\eta} < \nu \mathbb{L}^{gred}(\bm{y^{0,\eta}})$ holds. If it does, this indicates that $\theta^{\eta}$ is still relatively low, and we can safely add the relatively loose `greedy cut' to quickly raise $\theta$ without risking excluding the optimal solution that may obtained by the `augmented cut'. Otherwise, when $\theta^{\eta} > \nu \mathbb{L}^{gred}(\bm{y^{0,\eta}})$, it implies the current $\theta^\eta$ might be approaching the true value of $\mathbb{L}(\bm{y^{0,\eta}})$. In this case, we then resort to the more precise, yet (comparatively) computationally intensive, `augmented cut' by computing $\mathbb{L}^{aug}(\bm{y^{0,\eta}})$. Such an conditioned alternation between the two heuristics takes advantage of their relative strengthens, and therefore generally outperforms the cases where only one heuristic is consistently applied. This logic is implemented as Step 5 throughout the iterations in the final algorithmic framework outlined in \Cref{sec:cus_L_algo}.

\subsection{Warm Start for the L-Shaped Algorithm}\label{sec:warm_start_algo}

In addition to the closed loop of the algorithm, we design a warm start strategy to initialize it from two ``naive'' first-stage decisions, representing two extreme cases. The first decision is to not fulfill any order at the current cycle, but defer all requests to the future, i.e., $y_m^0 = 0, \forall m \in M$. On the contrary, the second decision is to undertake as many requests as possible in the current cycle within the service capacity, minimizing the postponement to the future. In particular, the second ``naive'' decision solves the following deterministic problem (\textbf{P4}), which myopically minimizes the total operational cost by ignoring the potential future request arrivals. A high penalty $c_m^{T+1}$ is applied to orders not fulfilled immediately.

(\textbf{P4}):
\begin{subequations}\label{P4}
	\begin{align}
            \min_{\bm{x^0, y^0, r^0} \in \mathds{B}; \bm{u^0, v^0} \in \mathds{R^+}}\ & \alpha \cdot \sum_{m \in M} \left(  c_m^0 \cdot y_m^0 + c_m^{T+1} \cdot y_m^{T+1}  \right) \nonumber \\
            + & \beta \cdot \sum_{(i,j)\in A} \left(d_{ij} \cdot\sum_{k \in K} x_{ij}^{0,k}\right) + \gamma \cdot \sum_{k \in K} (v_O^{0,k} - v_D^{0,k}) + \kappa \cdot \sum_{k \in K} r^{0,k} 
        \end{align}
	\begin{align}
		\text{s.t.} \quad &  (1) - (8) \nonumber \\
            & y_m^0 + y_m^{T+1} = 1 & \forall m \in M
	\end{align}
\end{subequations}

\noindent Then, corresponding to the two ``naive'' decisions, we generate `augmented cuts' and add them to the root node of the branch-and-bound tree (i.e., Step 1 in \Cref{cus_algorithm}). 

Further, in our problem, when the fulfillment of pending orders exceeds the service capacity of drones, there is limited flexibility to delay orders for cost savings. Consequently, if the above deterministic program outcomes indicate infeasibility of fulfilling all orders in the current cycle, we shift to a more streamlined approach. Instead of implementing the sequential cut-adding schema, we simply add the `greedy cut' through all iterations. This is achieved by omitting Step 5 in \Cref{cus_algorithm} while replacing $\mathbb{L}^{aug}(\bm{y}^0)$ with $\mathbb{L}^{gred}(\bm{y}^0)$ in Step 6.

\subsection{Customized L-shaped Algorithm}\label{sec:cus_L_algo}
Given the above elucidation, we now assemble the components into an algorithm following the standard framework of the integer L-shaped method:
\noindent
\begin{algorithm}
    \caption{Integer L-shaped method for dynamic order fulfillment with heuristic recourse}\label{cus_algorithm}
\end{algorithm}\vspace{-2em}
\begin{tcolorbox}[breakable]
    \begin{enumerate}\addtocounter{enumi}{-1}
        \item Set $\eta:=0$.
            \begin{itemize}
                \item Let $\bm{y}^{0,\eta}$ first equal to $0$ to compute $\mathbb{L}^{aug}(\bm{y^{0,\eta}})$ as the lower bound $L$ used for all subsequent optimality cuts. Next, solve (\textbf{P4}) to update $\bm{y^{0,\eta}}$, and initialize the best objective value $\bar{Z}=\text{Obj.}\ (\ref{P4}a) - \alpha \cdot \sum_{m \in M} c_m^{T+1} \cdot y_m^{T+1} + \mathbb{L}^{aug}(\bm{y^{0,\eta}})$.
                \item Create a single pendant node representing the master problem (\textbf{MP}) with two `augmented cuts' derived from the two ``naive'' solutions, and append it as the root node to the branch-and-bound tree.
            \end{itemize}
        \item Select a pendant node in the branching tree. If none exists, stop.
        \item Set $\eta:=\eta+1$; solve the current master problem (\textbf{MP}). 
        \begin{itemize}
            \item If the problem is infeasible, fathom the current node and loop back to Step 1.
            \item Else, let ($\bm{y^{0,\eta}}$, $\theta^\eta$) be the current optimal solution, and proceed to Step 3.
        \end{itemize}
        \item Compare the objective value obtained for (\textbf{MP}) with the best so far, i.e., $\text{Obj.}\ (\ref{MP}a)$ versus $\bar{z}$:
        \begin{itemize}
            \item If $\text{Obj.}\ (\ref{MP}a) \geq \bar{Z}$, there is no need to further branch on this node. Fathom the current problem and loop back to Step 1. 
            \item Else, proceed to Step 4.
        \end{itemize}
        \item Check the integrality of the current solution to (\textbf{MP}): \label{algo_step4}
        \begin{itemize}
            \item If any element of $\bm{y^{0,\eta}}$ violates the integrality constraint, generate two new branches following the standard branch-and-bound procedure. Append these two new nodes to the branching tree and loop back to Step 1.
            \item Else if all $\bm{y^{0,\eta}}$ are binary, resolve the master problem (\textbf{MP}) once again by enforcing both $\bm{x^{0,\eta}}$ and $\bm{r^{0,\eta}}$ to be binary.
            \begin{itemize}
                \item If the problem becomes infeasible, indicating that there is no attainable fulfillment plan, fathom this infeasible node, and loop back to Step 1.
                \item Else, it means that feasible plans exist under the current solution of $\bm{y^{0,\eta}}$, proceed to Step 5.
            \end{itemize}
        \end{itemize}
        \item Compute $\mathbb{L}^{gred}(\bm{y^{0,\eta}})$. 
        \begin{itemize}
            \item If $\theta^\eta < \nu \mathbb{L}^{gred}(\bm{y^{0,\eta}})$, then add a `greedy cut' and return to Step 2.
            \item Else, continue with Step 6.
        \end{itemize}
        \item Compute $\mathbb{L}^{aug}(\bm{y^{0,\eta}})$ and $Z^\eta = \text{Obj.}\ (\ref{MP}a) - \theta^\eta + \mathbb{L}^{aug}(\bm{y^{0,\eta}})$. If $Z^\eta < \bar{Z}$, update $\bar{Z} = Z^\eta$ and record current $\bm{y^{0,\eta}}$ as the best solution so far.
        \begin{itemize}
            \item If $\theta^\eta \geq \mathbb{L}^{aug}(\bm{y^{0,\eta}})$, then fathom the current node and loop back to Step 1.
            \item Else, append the `augmented cut' to the master problem and return to Step 2.
        \end{itemize}
    \end{enumerate}
\end{tcolorbox}

As a summary, in comparison to the standard integer L-shaped framework proposed by \cite{laporte1993integer}, \Cref{cus_algorithm} has introduced the following modifications to enhance the algorithmic efficiency:
\begin{itemize}

\item Instead of considering all first-stage binary variables, we build the branching tree and construct the optimality cut solely based on the fulfillment decision $\bm{y^0}$, excluding those related to the dispatching and routing of drones, which substantially slim down the tree size. Then, to ensure the feasibility of $\bm{y^0}$ when enforcing the integrality of other variables, an extra step of verification is included as Step \ref{algo_step4}.

\item A two-tier structure is introduced for the generation of optimality cuts. In the algorithm, we first add the `greedy cut' in Step 5, which maybe loose but can be computed quickly. After the loose cut improves the second-stage approximation $\theta$ to a certain level, we then gear-shift to the tighter `augmented cuts' in Step 6 for a finer approximation. The key idea is the computation of $\mathbb{L}^{gred}(\bm{y^{0,\eta}})$ is much faster than that of $\mathbb{L}^{aug}(\bm{y^{0,\eta}})$, allowing us to use the former as a rough filter to exclude some suboptimal solutions without incurring significant computational costs. 

\item A hyperparameter $\nu\in [0,1]$ is incorporated in Step 5 to act as the filter for deciding gear-switch for optimality cuts. A higher value of $\nu$ favors the `greedy cut' over the `augmented cut', thereby expediting the program's processing. However, a higher $\nu$ might violate the pre-condition set for \Cref{proposition_3}, potentially excluding the optimal solutions. According to Appendix \ref{sec:appd_subproblem_acc}, our experiments indicate that using only the `greedy cuts' can yield solutions with optimality gaps being less than 10\%. Concluding the trials, we adopt a moderate value setting $\nu=0.8$ in our numerical cases.
\end{itemize}

\noindent It is noteworthy that by setting $\nu=0$ and solving $\mathbb{L}(\bm{y^0})$ using the exact method, the proposed algorithm degenerates to the standard L-shaped method.

\section{Computational Performance Analysis}\label{sec:computational_analysis}
We undertake a comprehensive evaluation of the performance of our customized heuristic L-shaped method for solving the two-stage order fulfillment problem (\textbf{EF}). This evaluation is conducted in comparison with two benchmarks: the exact L-shaped method and the brute-force solution using Gurobi. Our primary focus centers on assessing both the quality of the solutions generated and the computational efficiency across these methodologies. Further, we perform a sensitivity analysis of our proposed two-tiered heuristic algorithm, examining its performance in relation to variations in the number of simulated scenarios and the volume of requests. All these methods are implemented in Python and deployed with Gurobi 9.5 on MacBook Pro with Apple M1 Pro Chip 3.22 GHz and 16 GB RAM. 

\subsection{Data Inputs} \label{data_generation}

To facilitate the computational experiments, a series of simulated data sets are created. Drones are parameterized according to their operating conditions under a free-flight mode. Drones navigate based on Euclidean distances, allowing for direct travel between points without the constraints of pre-defined pathways. The decision-making horizon spans three delivery cycles, encompassing the current cycle and two subsequent cycles to account for potential demand profiles. Requests not fulfilled within the three cycles are subject to high penalties in operational costs. Additionally, to provide an indication of the problem's size, we present the number of variables ($N_{vars}$) and constraints ($N_{cons}$). It should be noted that our framework is inherently scalable and can accommodate planning horizons of any length given. Table \ref{notation2} overviews the notations of parametric settings in simulations.

\begin{table}[ht!]
	\centering
	\caption{Notation list of simulated problem parameters}
	\label{notation2}
	\begin{tabular}{c l} 
	\hline
	Notation & Description \\
	\hline\hline
	$\left | M \right |$ & Number of first-stage requests\\
	$\left | N \right |$ & Number of second-stage requests\\
	$\left | K \right |$ & Number of drones\\
	$ T $ & Number of delivery cycles considered in the second stage\\	
	$\left |\Omega \right |$ & Number of simulated scenarios\\
	$N_{vars}$ & Number of variables for the SAA formulation of (\textbf{EF})\\
	$N_{cons}$ & Number of constraints for the SAA formulation of (\textbf{EF})\\
	\hline
	\end{tabular}
\end{table}

The experiments assume a square area of $2\times2$ square miles, with the service depot located at the center\footnote{The Federal Aviation Administration (FAA) permits a maximum speed of 100 mph for drones, enabling them to efficiently traverse inside a $2\times2$ square-mile area back and forth within the 10-minute delivery cycle, including takeoff and landing: \url{https://www.faa.gov/newsroom/small-unmanned-aircraft-systems-uas-regulations-part-107}.}. The x-y coordinates of future requests are simulated through uniform sampling within the [-1,2] mile range for each dimension. Requests falling into the negative coordinate regime are considered invalid and will be disregarded. This approach accommodates the randomness to the number of valid requests. In line with the specifications in \cite{jeong2019truck}, the load-dependent energy consumption function of drones $P(u_i)$ is assumed to be linear on the payload $u$, i.e., $P(u) = \lambda^p_0 + \lambda^p_1 \cdot u$, with $\lambda^p_0 = 0.1$ kW and $\lambda^p_1 = 0.1$ kW/kg. The weight of each package is uniformly generated within the range of $[5,10]$ kg. Additionally, the load capacity of each single drone is specified as $Q = 20$ kg, and its battery capacity is set at $E=0.5$ kW$\cdot$h \citep{xia2023branch}. The penalty for delayed order completion is represented by a piecewise function. For orders fulfilled within three cycles after their generation, the penalty is minimal, set at $c_m = 1$. However, for orders completed after $T+1$ cycles, the penalty escalates exponentially, following the formula $c_m = \lambda^c_c\cdot e^{\lambda^c_r\cdot \Delta t}$, where the coefficient $\lambda^c_c$ indicates the urgency level of each order, and $\lambda^c_r$ defines the sensitivity of the penalty to requests' elapsed time pending at the depot $\Delta t$.

When harmonizing the disparate sources of costs into the objective---spanning order fulfillment delays, vehicle depreciation, energy consumption, and labor supply---, it is imperative to ensure a coherent analysis in monetary terms (USD). To this end, we introduce specific weighting mechanisms for each cost category. For the penalty associated with order delays, we assign a value of $\alpha=\$0.05$/cycle. This translates to penalties of \$0.05 for regular orders completed within each incremental delivery cycle up to 30 minutes, then escalating significantly beyond this threshold to \$1.00, \$2.73, and \$7.42 for delays extending to 40, 50, and 60 minutes, respectively. Drawing on the cost analysis framework provided by Amazon Prime Air \citep{sudbury2016cost}, we assume an initial drone investment of \$5,000 with an operational lifespan of 500 hours, resulting in an amortized usage-based cost of $\beta=\$10$ per hour. For energy consumption costs, we set $\gamma=\$5$kW$\cdot$h for charging batteries. Additionally, there is a fixed labor cost of \$1 applied to each time of drone dispatch per delivery cycle.

\subsection{Performance Evaluations}

\begin{table}[b!]
    \footnotesize
    \centering
    \caption{Percent increases in the minimum costs (i.e., objective values) achieved across the five methods [Baseline: exact L-shaped method].\\ ($\left | M \right |=5$, $\left | N \right |=5$, $T=2$, $\left | K \right |=2$, $\left |\Omega \right |=10$, $N_{vars}=7,295$, $N_{cons}=12,896$).}\label{tab_comparison_obj}
\begin{tabular}{c|ccccc}
\hline
\multirow{2}{*}{Instance} & \multicolumn{5}{c}{Objective Value Gap} \\ \cline{2-6}
                          & \multicolumn{1}{c|}{Greedy Heu.} & \multicolumn{1}{c|}{Two-Tier Heu.} & \multicolumn{1}{c|}{Augmented Heu.} & \multicolumn{1}{c|}{Exact L-Shaped} & Brute-Force (MIP Gap) \\ \hline
1                         & \multicolumn{1}{c|}{3.42\%}    & \multicolumn{1}{c|}{0.47\%}          & \multicolumn{1}{c|}{0.47\%}    & \multicolumn{1}{c|}{0}              & 0.04\%  (57.50\%)                                                                             \\
2                         & \multicolumn{1}{c|}{1.29\%}    & \multicolumn{1}{c|}{0.38\%}          & \multicolumn{1}{c|}{0.38\%}    & \multicolumn{1}{c|}{0}              & 0.38\%  (55.74\%)                                                                             \\
3                         & \multicolumn{1}{c|}{2.07\%}    & \multicolumn{1}{c|}{0.07\%}          & \multicolumn{1}{c|}{0.07\%}    & \multicolumn{1}{c|}{0}              & 0.00\%  (55.55\%)                                                                             \\
4                         & \multicolumn{1}{c|}{2.39\%}    & \multicolumn{1}{c|}{0.60\%}          & \multicolumn{1}{c|}{0.60\%}    & \multicolumn{1}{c|}{0}              & 0.72\%  (57.99\%)                                                                             \\
5                         & \multicolumn{1}{c|}{1.19\%}    & \multicolumn{1}{c|}{0.12\%}          & \multicolumn{1}{c|}{0.12\%}    & \multicolumn{1}{c|}{0}              & 0.00\%  (59.39\%)                                                                             \\
6                         & \multicolumn{1}{c|}{0.44\%}    & \multicolumn{1}{c|}{0.12\%}          & \multicolumn{1}{c|}{0.12\%}    & \multicolumn{1}{c|}{0}              & 4.24\%  (60.27\%)                                                                             \\
7                         & \multicolumn{1}{c|}{5.13\%}    & \multicolumn{1}{c|}{0.70\%}          & \multicolumn{1}{c|}{0.70\%}    & \multicolumn{1}{c|}{0}              & 8.42\%  (58.58\%)                                                                             \\
8                         & \multicolumn{1}{c|}{1.98\%}    & \multicolumn{1}{c|}{0.68\%}          & \multicolumn{1}{c|}{0.68\%}    & \multicolumn{1}{c|}{0}              & 0.00\%  (58.09\%)                                                                             \\
9                         & \multicolumn{1}{c|}{2.85\%}    & \multicolumn{1}{c|}{0.94\%}          & \multicolumn{1}{c|}{0.94\%}    & \multicolumn{1}{c|}{0}              & 0.00\%  (58.66\%)                                                                             \\
10                        & \multicolumn{1}{c|}{2.08\%}    & \multicolumn{1}{c|}{0.07\%}          & \multicolumn{1}{c|}{0.07\%}    & \multicolumn{1}{c|}{0}              & 2.57\%  (58.41\%)                                                                             \\ \hline
Avg.                   & \multicolumn{1}{c|}{2.28\%}    & \multicolumn{1}{c|}{0.42\%}          & \multicolumn{1}{c|}{0.42\%}    & \multicolumn{1}{c|}{0}              & 1.64\%  (58.02\%)                                                                            \\
St. Dev.                   & \multicolumn{1}{c|}{1.31\%}    & \multicolumn{1}{c|}{0.31\%}          & \multicolumn{1}{c|}{0.31\%}    & \multicolumn{1}{c|}{0}              & 2.78\%  (1.47\%)                                                                              \\ \hline
\end{tabular}
\end{table}

\begin{table}[b!]
    \small
    \centering
    \caption{Running time across the five methods\\ ($\left | M \right |=5$, $\left | N \right |=5$, $T=2$, $\left | K \right |=2$, $\left |\Omega \right |=10$, $N_{vars}=7,295$, $N_{cons}=12,896$).}\label{tab_comparison_time}
\begin{tabular}{c|ccccc}
\hline
\multirow{2}{*}{Instance} & \multicolumn{5}{c}{Running Time (seconds)}                                                                                                                           \\ \cline{2-6} 
                          & \multicolumn{1}{c|}{Greedy Heu.} & \multicolumn{1}{c|}{Two-Tier Heu.} & \multicolumn{1}{c|}{Augmented Heu.} & \multicolumn{1}{c|}{Exact L-Shaped} & Brute-Force \\ \hline
1                         & \multicolumn{1}{c|}{8}         & \multicolumn{1}{c|}{53}              & \multicolumn{1}{c|}{71}        & \multicolumn{1}{c|}{1,199}           & \textgreater{}3,600    \\
2                         & \multicolumn{1}{c|}{8}         & \multicolumn{1}{c|}{34}              & \multicolumn{1}{c|}{69}        & \multicolumn{1}{c|}{1,355}           & \textgreater{}3,600    \\
3                         & \multicolumn{1}{c|}{8}         & \multicolumn{1}{c|}{60}              & \multicolumn{1}{c|}{71}        & \multicolumn{1}{c|}{1,168}           & \textgreater{}3,600    \\
4                         & \multicolumn{1}{c|}{8}         & \multicolumn{1}{c|}{62}              & \multicolumn{1}{c|}{71}        & \multicolumn{1}{c|}{1,130}           & \textgreater{}3,600    \\
5                         & \multicolumn{1}{c|}{8}         & \multicolumn{1}{c|}{62}              & \multicolumn{1}{c|}{71}        & \multicolumn{1}{c|}{1,433}           & \textgreater{}3,600    \\
6                         & \multicolumn{1}{c|}{8}         & \multicolumn{1}{c|}{57}              & \multicolumn{1}{c|}{72}        & \multicolumn{1}{c|}{1,566}           & \textgreater{}3,600    \\
7                         & \multicolumn{1}{c|}{8}         & \multicolumn{1}{c|}{53}              & \multicolumn{1}{c|}{72}        & \multicolumn{1}{c|}{1,517}           & \textgreater{}3,600    \\
8                         & \multicolumn{1}{c|}{8}         & \multicolumn{1}{c|}{61}              & \multicolumn{1}{c|}{71}        & \multicolumn{1}{c|}{1,214}           & \textgreater{}3,600    \\
9                         & \multicolumn{1}{c|}{8}         & \multicolumn{1}{c|}{60}              & \multicolumn{1}{c|}{71}        & \multicolumn{1}{c|}{1,215}           & \textgreater{}3,600    \\
10                        & \multicolumn{1}{c|}{8}         & \multicolumn{1}{c|}{54}              & \multicolumn{1}{c|}{71}        & \multicolumn{1}{c|}{1,072}           & \textgreater{}3600    \\ \hline
Avg.                   & \multicolumn{1}{c|}{8}         & \multicolumn{1}{c|}{56}              & \multicolumn{1}{c|}{71}        & \multicolumn{1}{c|}{1,287}           & \textgreater{}3,600    \\
St. Dev.                   & \multicolumn{1}{c|}{0}         & \multicolumn{1}{c|}{8}               & \multicolumn{1}{c|}{1}         & \multicolumn{1}{c|}{170}            & -                     \\ \hline
\end{tabular}
\end{table}

We simulate ten instances with the following setup, $\left | M \right |=5$, $\left | N \right |=5$, $T=2$, $\left | K \right |=2$, and $\left |\Omega \right |=10$, and solve the order dispatching problem using five different approaches. The quality of solutions and computational efficiency resulting from each method are reported and compared in \Cref{tab_comparison_obj,tab_comparison_time}, respectively. The `greedy heuristic', `augmented heuristic', and `exact L-shaped' methods are implemented using the customized integer L-shaped algorithm described in \Cref{sec:cus_L_algo}, with Step 5 omitted and then $\mathbb{L}^{gred}(\bm{y^0})$, $\mathbb{L}^{aug}(\bm{y^0})$, and the exact $\mathbb{L}(\bm{y^0})$ used in respective cases. The `two-tier heuristic' method strictly executes \Cref{cus_algorithm}, serving as an intermediate between the two heuristic methods. Additionally, the `brute-force' method directly solves the extended form (\textbf{EF}) using Gurobi 9.5, with computational time capped at one hour. If the solution time exceeds this limit, we terminate the program and report the immediate outcomes. We use the objective value obtained from the `exact L-shaped' method as the benchmark, given its guarantee of global optimal solutions. Then, by comparing the objective values achieved by each method against this benchmark, we calculate the metric `objective value gap'. Particularly, for the `brute-force' method, the metric `MIP gap' reported by Gurobi is also included for reference. This metric reflects the quality of the solution that Gurobi achieves upon termination. A smaller gap indicates a solution closer to optimality, with the optimal solution being indicated by a zero gap. 

It is evident that solving the SAA of the two-stage stochastic program directly using Gurobi is impractical. With 7,295 decision variables and 12,896 constraints, Gurobi constantly experiences difficulty in finding the optimal solution within the predefined one-hour time cap. Since adding one extra scenario roughly introduces hundreds of additional variables, the computational time can easily explode when using Gurobi for instances with a large number of scenarios. Moreover, while the MIP gaps hover consistently around 58\% upon termination, the objective values exhibit considerable variance in distances from the optimal solution, indicating unstable solution quality. In contrast, all four customized L-shaped methods drastically reduce the running time to a manageable level. Furthermore, the three heuristic methods enhance computational efficiency by over 20 times compared to the exact customized L-shaped method. As expected, the `greedy heuristic' method emerges as the quickest in solving the problem, albeit incurring a relatively higher objective gap. The `two-tier heuristic', which alternates between `greedy cut' and `augmented cut', matches the same solution quality as the `augmented heuristic' but with a 20\% improvement in running time. Later, we will show more examples that illustrate the superior performance of the L-shaped methods in handling realistic problems. Given Gurobi's inability to achieve an optimal solution within a reasonable time frame for those cases, it is excluded from further comparisons.

\subsection{Sensitivity to the Number of Scenarios}

\begin{figure}[t!]
	\centering
	\begin{subfigure}[b]{.5\textwidth}
		\centering
		\includegraphics[width=1\linewidth,height=0.3\textheight]{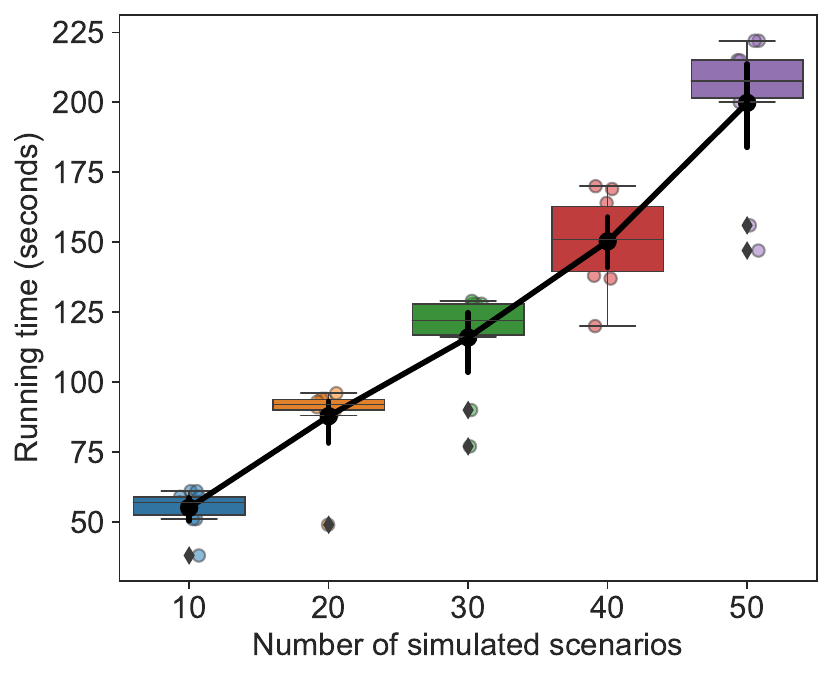}
		\caption{Running Time (seconds)}
	\end{subfigure}%
	\begin{subfigure}[b]{.5\textwidth}
		\centering
		\includegraphics[width=1\linewidth,height=0.3\textheight]{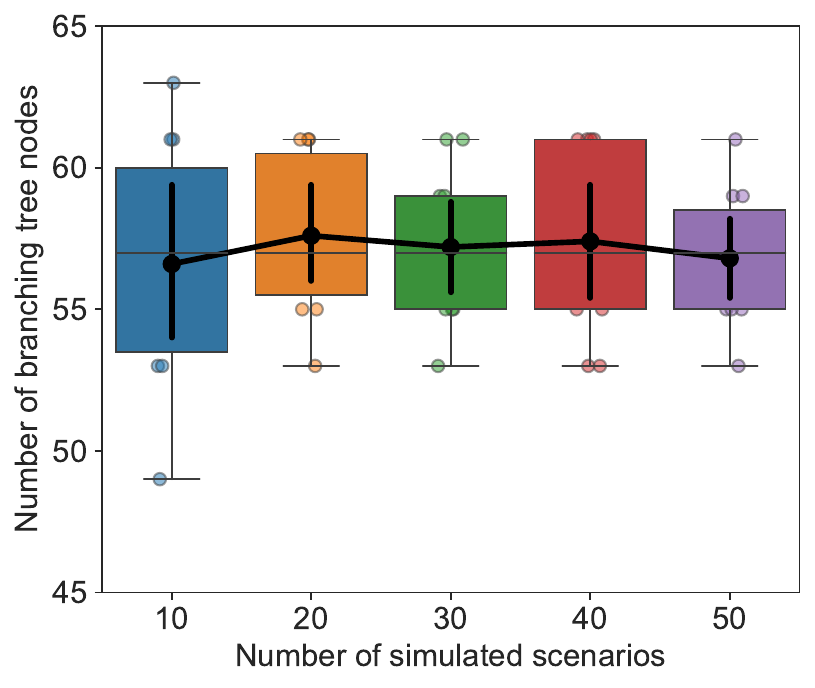}
		\caption{Number of Branching Tree Nodes}
	\end{subfigure}
	\caption{Changes in (a) running time (seconds) and (b) the number of branching tree nodes with respect to the number of simulated scenarios}
	\label{fig_n_scen}
\end{figure}

We proceed with a sensitivity analysis regarding the number of simulated scenarios. All experiments for this analysis are handled by the `two-tier heuristic' method. Figure \ref{fig_n_scen} shows the box plots of the running time and the number of branching tree nodes for problems of varying scenario sizes. Each box plot is derived from ten replications. The number of branching tree nodes refers to the total number of nodes added to the branching tree within the integer L-shaped method. Each time of branching, based on a non-integer variable, results in two additional nodes appended to the branching tree. In general, a larger branching tree size corresponds to more iterations to execute, while the computing time of each iteration is dictated by the time needed for the construction of optimality cuts on each node. As a result, the total solution time for each instance depends on both the size of the branching tree and the complexity of constructing optimality cuts. Therefore, both running time and the number of branching nodes are reported as indicators of computational efficiency.

Figure \ref{fig_n_scen} illustrates that the running time of the customized integer L-shaped method scales linearly with the growing number of scenarios, while the size of branching tree nodes exhibits minimal variation. This phenomenon can be attributed to the fact that the branching tree size is predominantly influenced by the number of integer variables in the first stage, unaffected by the number of scenarios in the second stage. Therefore, the augmentation in scenario count primarily affects the efficiency of computing optimality cuts. Specifically, as evidenced by \Cref{fig_n_scen}a, the time dedicated to the computation of optimality cuts is linearly proportional to the number of scenarios. In fact, if the number of computer cores exceeds that of scenarios, individual realization of the subproblem $\mathbb{Q}(\bm{y^{0,\eta}}, \omega)$ can be processed entirely in parallel. In such cases, the running time would remain constant for the varying numbers of simulated scenarios. 

\subsection{Sensitivity to the Number of Requests}

\begin{figure}[b!]
	\centering
	\begin{subfigure}[b]{.5\textwidth}
		\centering
		\includegraphics[width=0.9\linewidth,height=0.3\textheight]{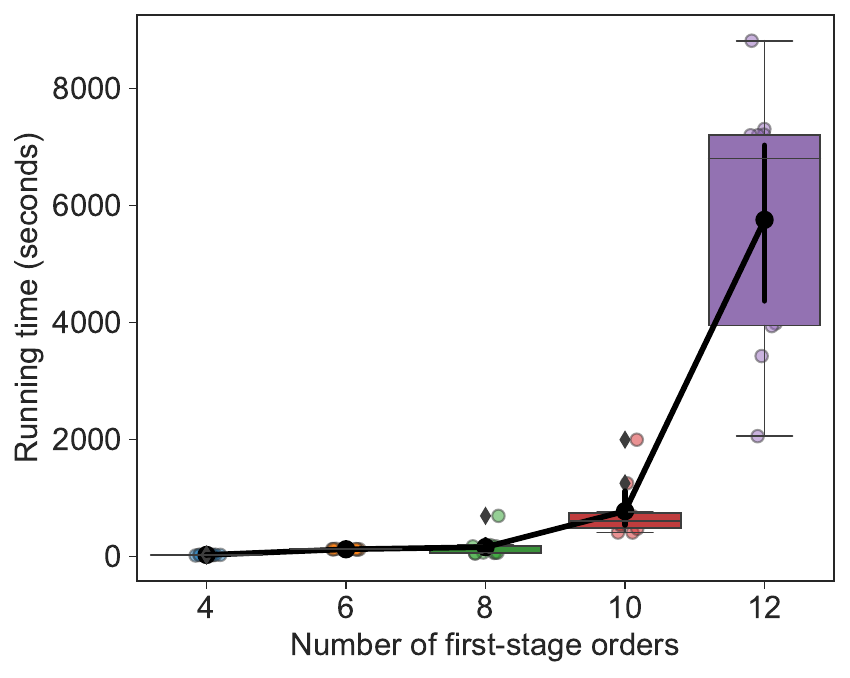}
		\caption{Running Time (seconds)}
	\end{subfigure}%
	\begin{subfigure}[b]{.5\textwidth}
		\centering
		\includegraphics[width=0.9\linewidth,height=0.3\textheight]{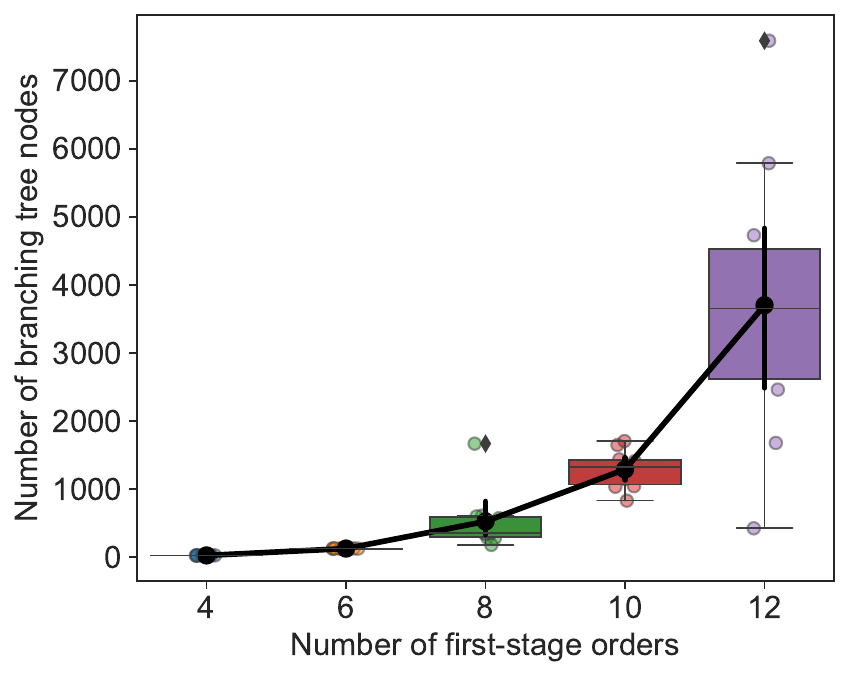}
		\caption{Number of Branching Tree Nodes}
	\end{subfigure}
	\caption{Changes in (a) running time (seconds) and (b) the number of branching tree nodes with respect to the number of orders $\left | M \right |$ considered in the first stage.}
	\label{boxplot_M}
\end{figure}

\begin{figure}[b!]
	\centering
	\begin{subfigure}[b]{.5\textwidth}
		\centering
		\includegraphics[width=0.9\linewidth,height=0.3\textheight]{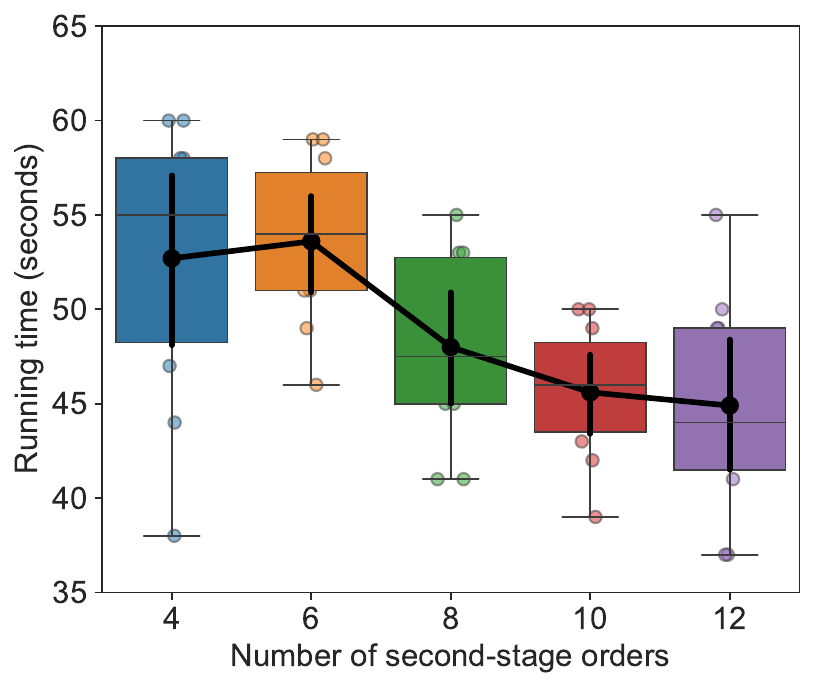}
		\caption{Running Time (seconds)}
	\end{subfigure}%
	\begin{subfigure}[b]{.5\textwidth}
		\centering
		\includegraphics[width=0.9\linewidth,height=0.3\textheight]{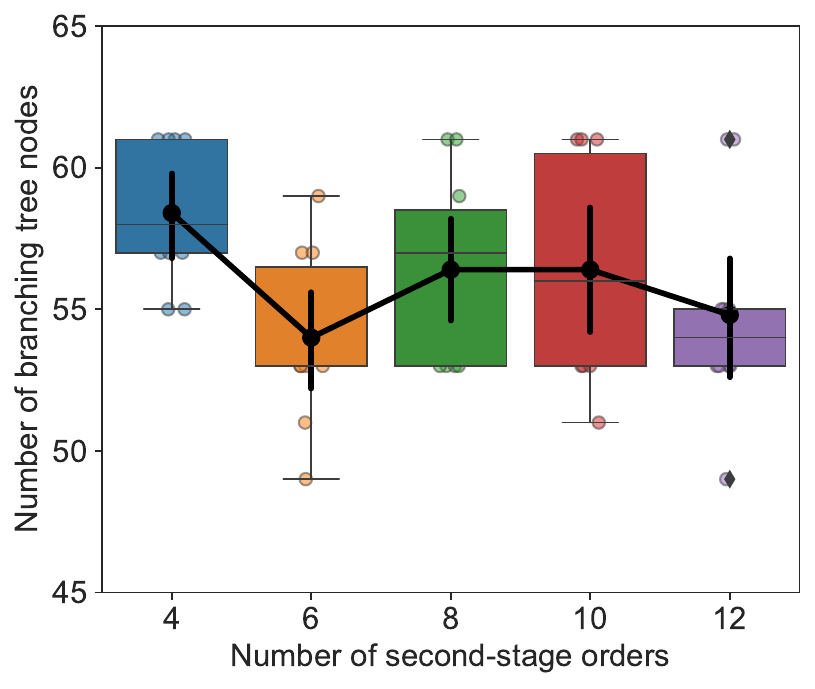}
		\caption{Number of Branching Tree Nodes}
	\end{subfigure}
	\caption{Changes in (a) running time (seconds) and (b) the number of branching tree nodes with respect to the number of orders $\left | N \right |$ in the second stage.}
	\label{boxplot_N}
\end{figure}

Another dimension that may influence the performance of the `two-tier heuristic' method is the volume of order requests. Figures \ref{boxplot_M} and \ref{boxplot_N} present the responsiveness of our customized L-shaped method's running time and branching tree size to changes in the number of requests in the current delivery wave $\left | M \right |$ and future waves $\left | N \right |$, respectively.

Figure \ref{boxplot_M} first reveals that both the running time and the number of branching tree nodes exhibit an approximate exponential increase with the volume of requests in the first stage. Given that the branching algorithm applies to the first-stage variable $\bm{y^0}$, such exponential responses are predictable. In particular, a transition occurs when the delivery volume exceeds the service capacity of delivery drones, prompting the algorithm to shift from the `two-tier heuristic' to using the `greedy heuristic' exclusively. For example, in our experiments, when the cumulative weight of $|M|=8$ packages exceeds the carrying capacity of drones, this transition is triggered, resulting in the slight decline in running time compared to the instances with $|M|=6$. That being said, as we further increase $|M|$ to 10 orders, the algorithm faces significant computational challenges in handling the growing VRP. As discussed in \Cref{sec:introduction}, drone delivery is particularly suitable for scenarios characterized by sparse and uncertain demand. However, when order density increases, the limited capacity of drones necessitates frequent round trips between the depot and customers. In such cases, the economies of scale associated with traditional truck delivery become more favorable. Meanwhile, to address large-scale VRPs, a common approach involves segmenting the problem into smaller, more manageable segments. For example, orders can be clustered in a preprocessing stage, and then these clusters are assigned to drones, resulting in a more compact, cluster-based MILP problem formulation \citep{dondo2007cluster}. Another widespread method partitions the area into regions, addressing the VRP for each region independently \citep{taillard1993parallel, reimann2004d}. While these methods are beyond the scope of this paper, they can be considered for future extensions.

Comparatively, Figure \ref{boxplot_N} shows that the running time remains relatively constant concerning the volume of potential orders $\left | N \right |$. As explained in Appendix \ref{sec:appd_subproblem_acc}, the increase in orders within the subproblem barely affects the running time of both the `greedy heuristic' and `augmented heuristic'. Since the cost of computing one optimality cut and the size of the branching tree remain stable, increasing number of requests in the second stage does not adversely affect the efficiency of our customized L-shaped method. 

In summary, adding one additional order to either the first-stage $\left | M \right |$ or the second-stage $\left | N \right |$ would both introduce thousands more integer variables to the problem. However, increasing the volume of potential orders would not affect the size of the branching tree, whilst a rise in the volume of existing orders can lead to a significant increase in the number of branching tree nodes to be computed. Consequently, the running time of our customized method exhibits an exponential increase with respect to $\left | M \right |$, but remains constant when varying $\left | N \right |$.

\section{Operational Performance Analysis}\label{sec:operational_analysis}

This section transitions our focus from computational performance to the operational performance established by the two-stage order fulfillment model. Another series of numerical experiments is conducted to further evaluate the effectiveness of our proposed stochastic model in enhancing operational efficiency in dynamic order fulfillment, and examine influential factors related to the various aspects of demand profiles. Specifically, the comparisons are made between the two-stage program, whose second stage accounts for potential order arrivals, and the myopic single-stage counterpart (\textbf{P4}) (i.e., with only the first stage). The discrepancies in their performances imply the values of considering future demand in operations, despite its uncertainty to varying extents. Given the superior computational performance of the customized L-shaped algorithm, as demonstrated by \Cref{sec:computational_analysis}, it is employed to conduct the numerical experiments in this section. The exact method is not used due to its prohibitively long computational time.

\subsection{Experimental Setup}

The ``ground truth'' arrivals of delivery requests is simulated over an extended period to experiment with the order fulfillment models on a cyclically rolling basis. The number of realized orders for each cycle is generated from a gamma distribution $\Gamma(\mu_\gamma, \sigma_\gamma)$ and rounded to an integer. The location of each order is randomly generated from a truncated normal distribution $N(\mu_n, \sigma_n^2)$ within a $2\times2$ square-mile map. Specifically, $\mu_\gamma,\mu_n$ and $\sigma_\gamma,\sigma_n$ are the means and standard deviations of the gamma and normal distributions, respectively. The package weights are still generated from a uniform distribution $U(a_u, b_u)$, where $a_u$ and $b_u$ indicate the minimum and maximum weights. The remaining parameters follow their values specified in \Cref{data_generation}.

In the context of the two-stage stochastic program, the depot manager possesses precise information regarding the existing orders and has knowledge of the probability distributions for demand profiles in subsequent cycles, i.e., knowing the statistics $\{\mu_\gamma, \sigma_\gamma, \mu_n, \sigma_n, a_u, b_u\}$ from historical service records. At each decision epoch, a set of $|\Omega|$ scenarios is generated based on the upcoming demand profile, feeding into the second-stage sub-problem to simulate future demand realizations. From cycle to cycle, we implement the two-stage program and then execute the first-stage solutions on the ``ground truth'' dataset; finally, we assess the resulting performance in comparison to that of the single-stage model. 

\subsection{Cost Saving Mechanisms}

\begin{figure}[b!]
	\centering
	\begin{subfigure}[b]{.4\textwidth}
		\centering
		\includegraphics[width=1\linewidth,height=0.24\textheight]{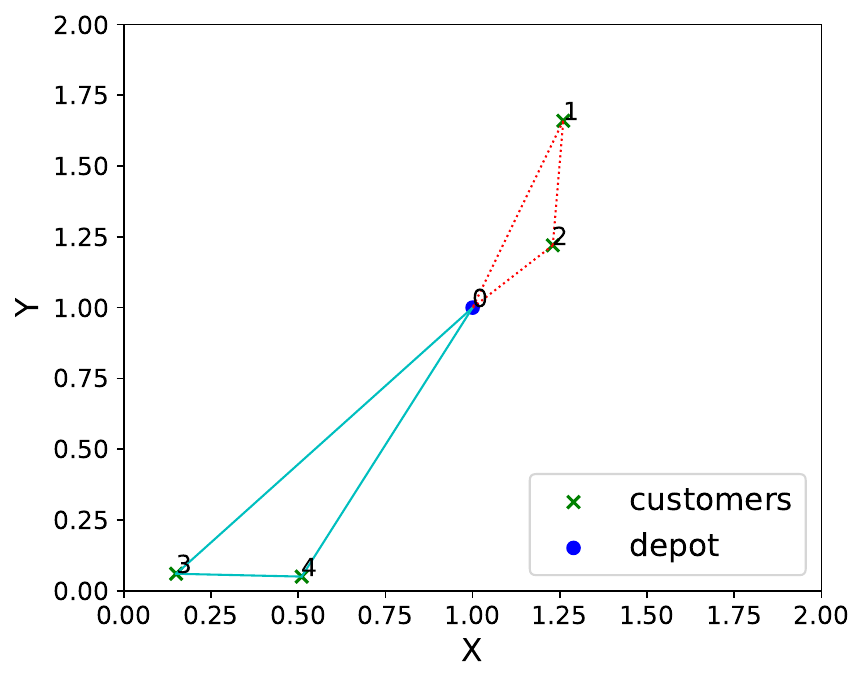}
		{\small Single-Stage / Cycle 1}
	\end{subfigure}%
	\begin{subfigure}[b]{.4\textwidth}
		\centering
		\includegraphics[width=1\linewidth,height=0.24\textheight]{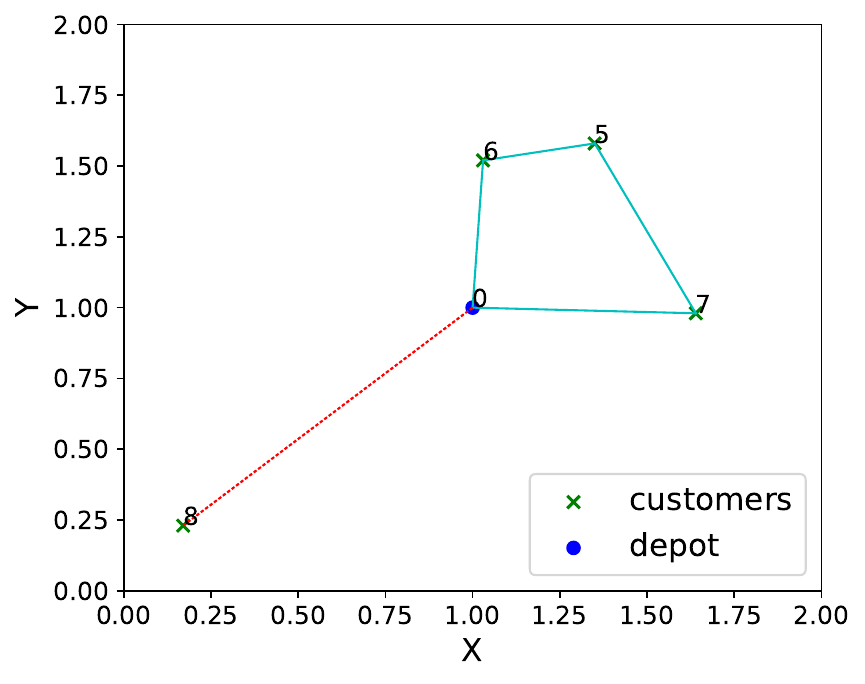}
		{\small Single-Stage / Cycle 2}
	\end{subfigure}\\
	\begin{subfigure}[b]{.4\textwidth}
		\centering
		\includegraphics[width=1\linewidth,height=0.24\textheight]{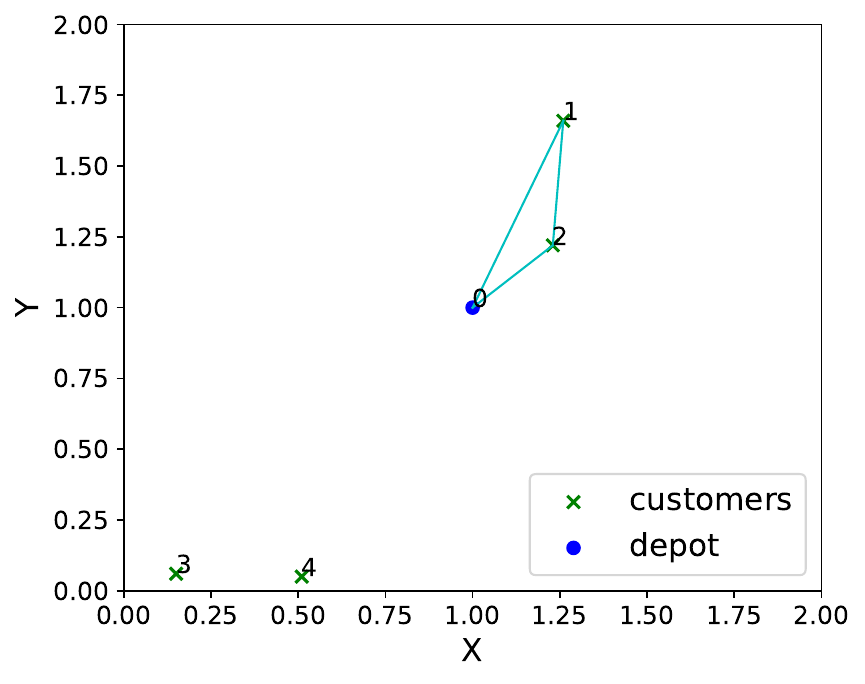}
		{\small Two-Stage / Cycle 1}
	\end{subfigure}%
	\begin{subfigure}[b]{.4\textwidth}
		\centering
		\includegraphics[width=1\linewidth,height=0.24\textheight]{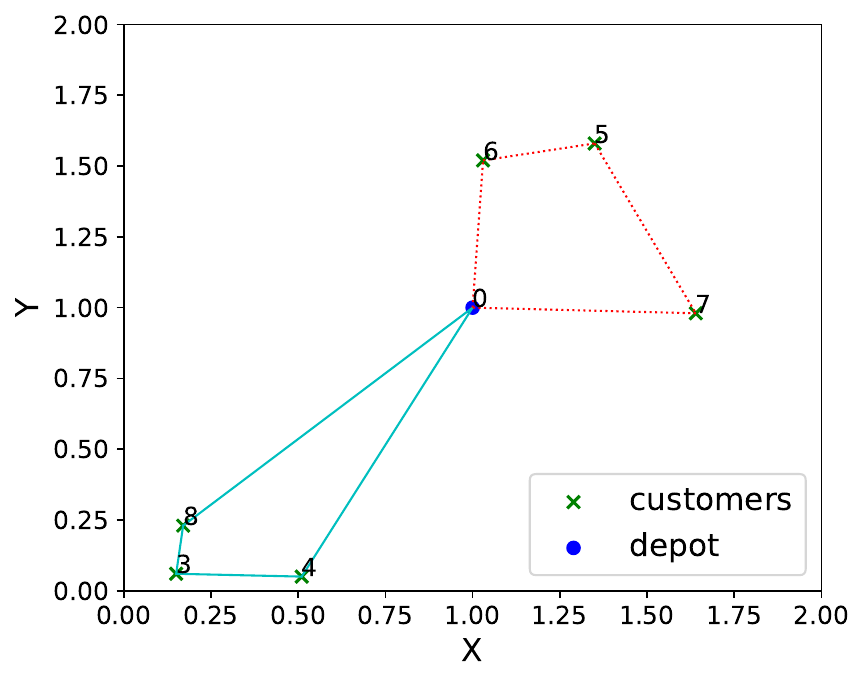}
		{\small Two-Stage / Cycle 2}
	\end{subfigure}
	\caption{An example showcasing trip reduction achieved by the two-stage program compared to the single-stage program. Node `0' represents the central depot, while other nodes depict access points for delivery. Solid and dotted arcs illustrate the optimal routes of individual drones. The result comparisons are presented as three-element tuples, in which the costs for the single-stage program, the two-stage program, and their differences are spelled out sequentially: \textit{Total cost - $(\$7.18,\ \$5.63\ | -21.6\%)$; Order delay penalty - $(\$0.40,\ \$0.40\ | +0.0\%)$; Distance-related cost - $(\$2.15,\ \$1.59\ | -26.0\%)$; Energy-related cost - $(\$0.64,\ \$0.64\ | -0.0\%)$; Drone-dispatching cost - $(\$4.00,\ \$3.00 |\ -25.0\%)$}. The absolute cost values duplicate the corresponding cost terms in the objective of the single- and two-stage programs in (\textbf{P4}) and (\textbf{MP}), respectively.
 }
	\label{fig_trip_sav}
\end{figure}

In general, the two-stage model integrates future delivery request information, thereby enhancing coordination in order fulfillment across consecutive cycles, allowing for collective and more efficient service. Before delving into quantitative analysis, we illustrate two specific mechanisms through which the two-stage model outperforms the single-stage counterpart in lowering operational costs.

\begin{figure}[b!]
	\centering
 	\begin{subfigure}[b]{.4\textwidth}
		\centering
		\includegraphics[width=1\linewidth,height=0.24\textheight]{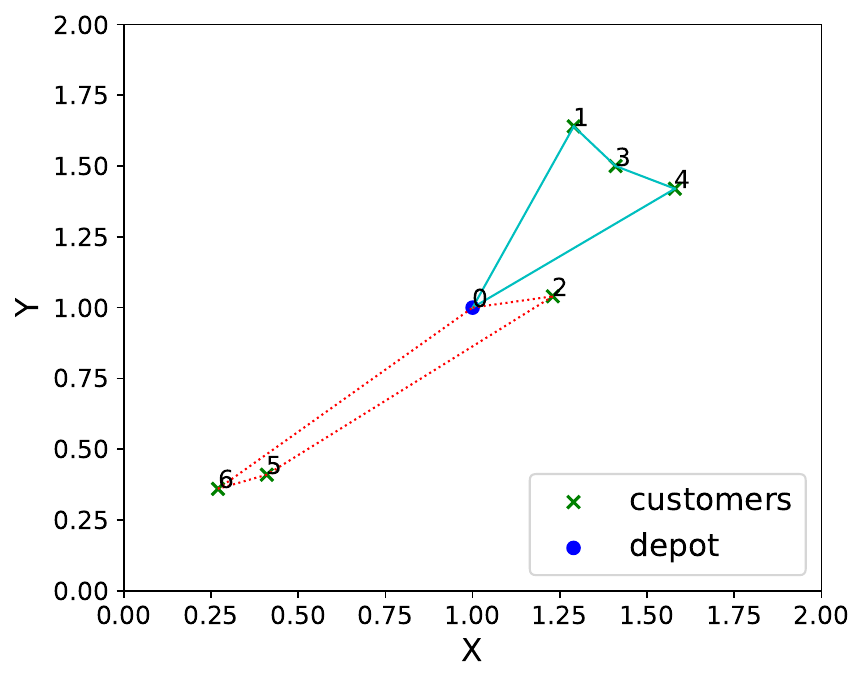}
		{\small Single-Stage / Cycle 1}
	\end{subfigure}%
	\begin{subfigure}[b]{.4\textwidth}
		\centering
		\includegraphics[width=1\linewidth,height=0.24\textheight]{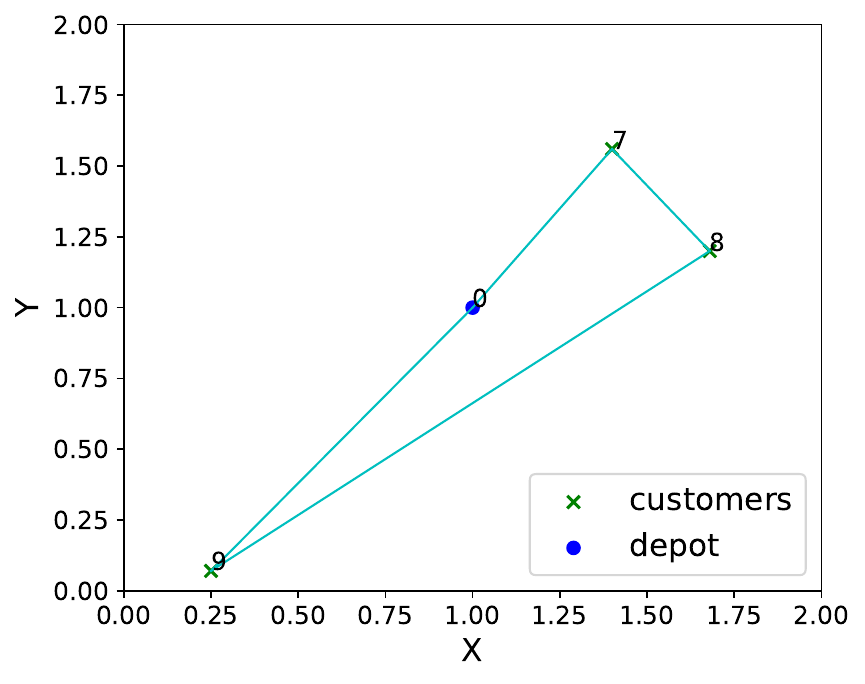}
		{\small Single-Stage / Cycle 2}
	\end{subfigure}\\
	\begin{subfigure}[b]{.4\textwidth}
		\centering
		\includegraphics[width=1\linewidth,height=0.24\textheight]{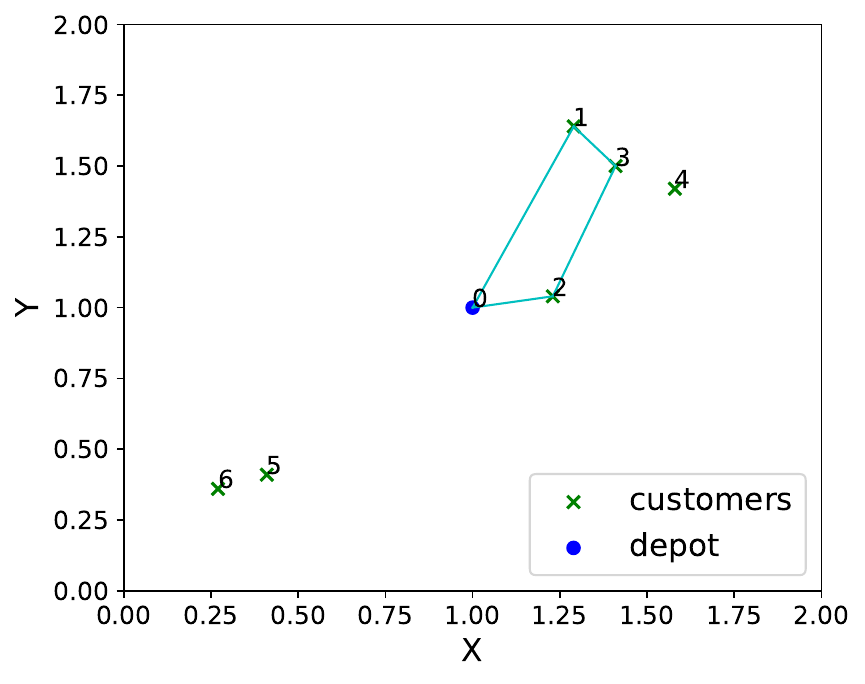}
		{\small Two-Stage / Cycle 1}
	\end{subfigure}%
	\begin{subfigure}[b]{.4\textwidth}
		\centering
		\includegraphics[width=1\linewidth,height=0.24\textheight]{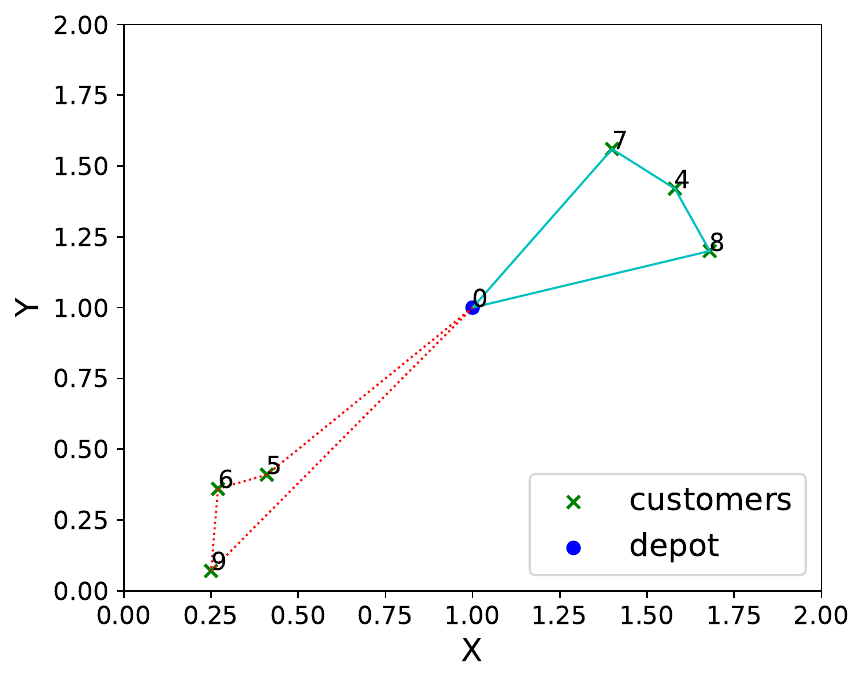}
		{\small Two-Stage / Cycle 2}
	\end{subfigure}
	\caption{An example of trip savings achieved by the stochastic program. The itemized costs resulting from the single- and two-stage programs, along with their differences, are again presented as three-element tuples: \textit{Total cost - $(\$6.34,\ \$5.53 \ | -12.8\%)$; Order delay penalty - $(\$0.45,\ \$0.45\ | +0.0\%)$; Distance-related cost - $(\$2.08,\ \$1.49\ | -28.7\%)$; Energy-related cost - $(\$0.81,\ \$0.59\ | -26.7\%)$; Drone-dispatching cost - $(\$3.00,\ \$3.00\ |\ -0.0\%)$}.
    }
	\label{fig_distance_sav}
\end{figure}

The first mechanism operates through \textit{trip reductions}, which, when applicable, can usually bring in significant cost savings. In cases where the two-stage stochastic program identifies a high potential of orders appearing near existing ones and anticipates that drones can handle them collectively, it defers these orders to the next cycle for simultaneous fulfillment. This strategic delay slightly heightens the risk of delay penalties in trade of substantial savings by reducing the number of delivery trips. \Cref{fig_trip_sav} illustrates a depot equipped with two identical drones practicing the trip-reduction mechanism. The solid and dotted trajectories depict the delivery tours of the two drones in individual cycles. In this instance, the single-stage model exhausts capacity by accepting as many existing orders as possible in each of the two cycles, resulting in a separate delivery leg to accommodate new arrivals of Request 8 in Cycle 2. In contrast, the two-stage model strategically delays Requests 3 and 4, fulfilling them alongside Request 8. This saves one trip, leading to a total saving of 21.6\%, with a 26.0\% reduction in traveling distance costs and a 25.0\% reduction in drone dispatch fixed costs, respectively.

But there may not always be opportunities to batch spatially close orders across different cycles. Besides trip reductions, the two-stage program can also optimize the arrangement for delivery tours, adhering to the second mechanism of \textit{trip savings}. Even with the same number of trips, strategically delaying orders enables distinct grouping and routing for order deliveries, creating opportunities for savings in distance and energy consumption. This mechanism is illustrated through \Cref{fig_distance_sav}. Unlike the previous example, both models suggest three delivery trips, indicating no opportunity for trip reductions. Nevertheless, the two-stage program still achieves a 12.8\% reduction in total costs. The single-stage program, while attempting to fulfill as many orders as possible, instructs one drone traversing the entire map in Cycle 2 to accommodate the remaining delivery orders. In contrast, the two-stage model intelligently batches a close group of orders in Cycle 1, avoiding extensive long-distance deliveries, and then fulfills the remaining orders in the second cycle, similar to Cycle 1 under the single-stage counterpart. Such a ``rearragneemt'' of fulfillment sequences reduces the cost of traveling distance and energy consumption by 28.7\% and 26.7\%, respectively. 

\subsection{Numerical Results and Analyses}

The aforementioned mechanisms of trip reductions and trip savings, showcased over two cycles, contribute to the superiority of the two-stage model in long-term operations. To illustrate, we conduct a series of experiments based on a realistic setting and adjust various parametric inputs to examine how the proposed two-stage order-dispatch model could perform in operations.

The experiments assume a two-hour lunch peak with 12 delivery cycles. To replicate spatial anisotropy, we configure the $2\times2$ square-mile map to include two clustered areas for delivery destinations surrounding the central depot. These clusters are respectively centered at coordinates $(\mu_n^{X1}, \mu_n^{Y1}) = (1.5, 1.5)$ miles and $(\mu_n^{X2}, \mu_n^{Y2})= (0.2, 0.2)$ miles. The average numbers of delivery requests for the two demand clusters across the 12 cycles are set to $\mu_\gamma^1 = \{ 1,2,3,3,4,3,3,4,3,3,2,1 \}$ and $\mu_\gamma^2 = \{ 1,1,1,2,2,1,1,2,2,1,1,1 \}$ to replicate the peaking and unpeaking process of service peak hours. The above set of parameteric values remains constant throughout the experiments, while the variances for these inputs are adjusted to explore different dimensions and levels of uncertainties. For each set of parameters, we conduct 10 replications of the experiment to provide an indicative assessment.

\begin{figure}[t!]
	\centering
        \begin{subfigure}[b]{.9\textwidth}
		\centering
		\includegraphics[width=1.0\textwidth]{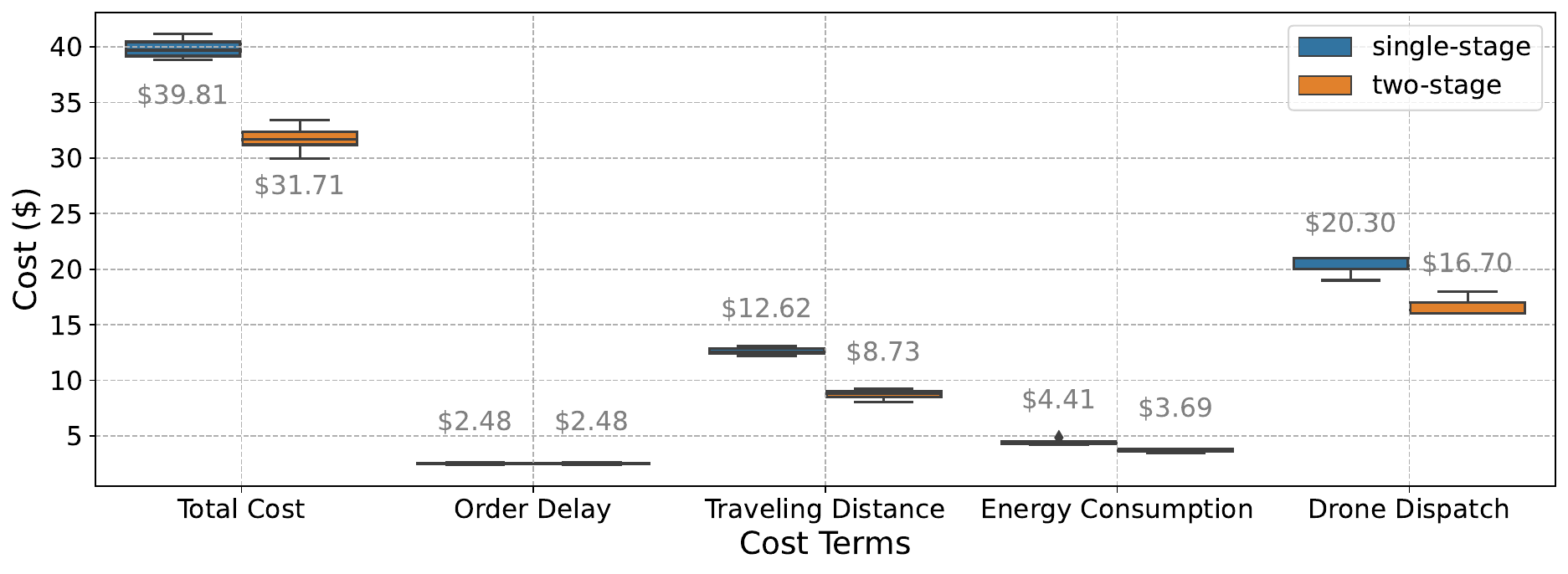}
    	\caption{Itemized costs (single-stage versus two-stage)}
            \label{Cost_saving_value}
	\end{subfigure}\\
	\begin{subfigure}[b]{.9\textwidth}
		\centering
		\includegraphics[width=1.0\textwidth]{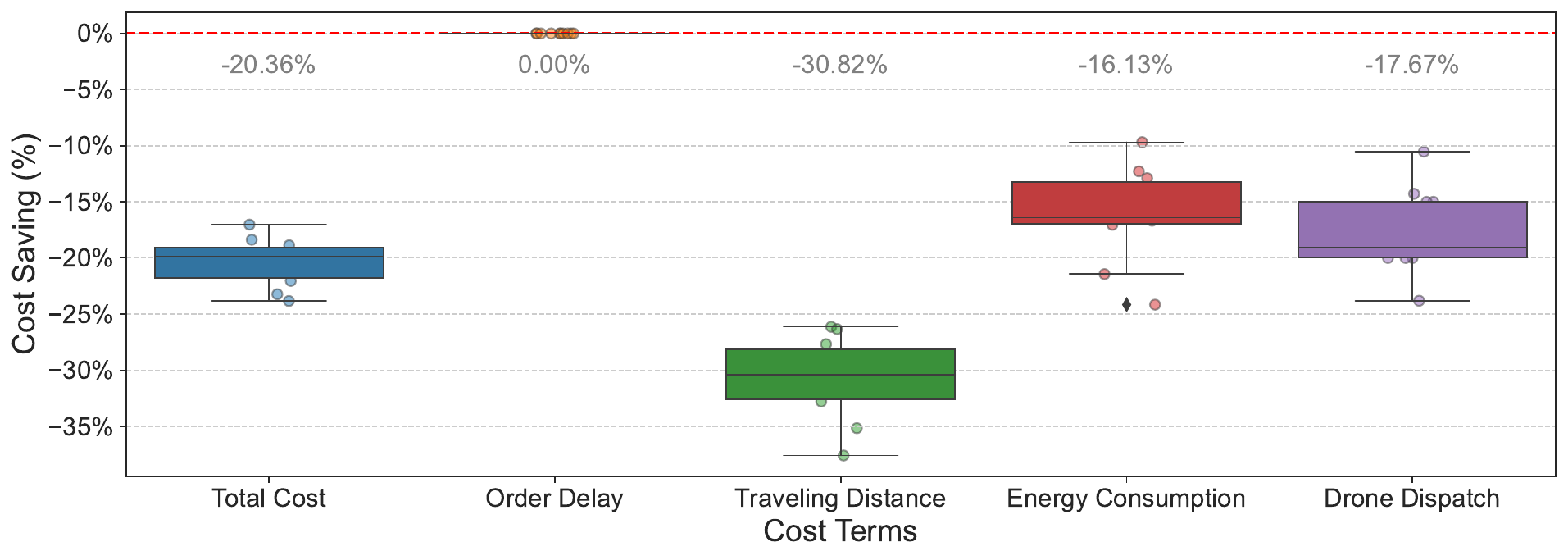}
    	\caption{Percentage cost reductions (single-stage $\Rightarrow$ two-stage)}
    	\label{Cost_saving_diff}
	\end{subfigure}
        \caption{Cost comparisons for two-hour operations with two-stage versus single-stage models (10 replications).}
        \label{fig_lunch peak_base}
\end{figure}

We begin by analyzing the results of a base scenario, where we set the standard deviation of the number of orders per cycle in each cluster as $\sigma_\gamma^1=\sigma_\gamma^2=0.1$, the standard derivation of the destination concentration for each cluster as $\sigma_n^1=\sigma_n^2=0.1$ mile, and the lower and upper bounds of the package weights as $a_u=5$ and $b_u=6$ kilogram. Intentionally, the uncertainty is kept relatively low for the base scenario, with a comprehensive sensitivity analysis for increased levels of uncertainty following subsequently. The results from implementing single-stage and two-stage models in the base scenario are summarized in \Cref{fig_lunch peak_base}. Figure \ref{Cost_saving_value} shows that in a low-uncertainty environment, the outcomes of both models are relatively stable. All the orders in this scenario are fulfilled within $T$ cycles without incurring additional delay penalty. But consistently, the two-stage model results in savings across various aspects of operational costs. As clearly presented in Figure \ref{Cost_saving_diff}, there is a substantial 20.36\% reduction in the total cost by switching from single- to two-stage models, accompanied by average reductions of 30.82\% in travel distance, 16.13\% in energy consumption, and 17.67\% in the frequency of drone dispatch. It is worth noting that among the four cost components, the frequency of drone dispatch is associated with the mechanism of trip reduction, while the travel distance and energy consumption are a joint reflection of trip reduction and trip saving mechanisms. The reduction in energy consumption is moderate due to a trade-off between the number of delivery trips and the package loads per trip. While shorter travel distances reduce energy use, fewer trips result in heavier loads per trip, which increases energy consumption.

\begin{figure}[b!]
	\centering
        \begin{subfigure}[b]{.9\textwidth}
            \centering
            \includegraphics[width=1.0\textwidth]{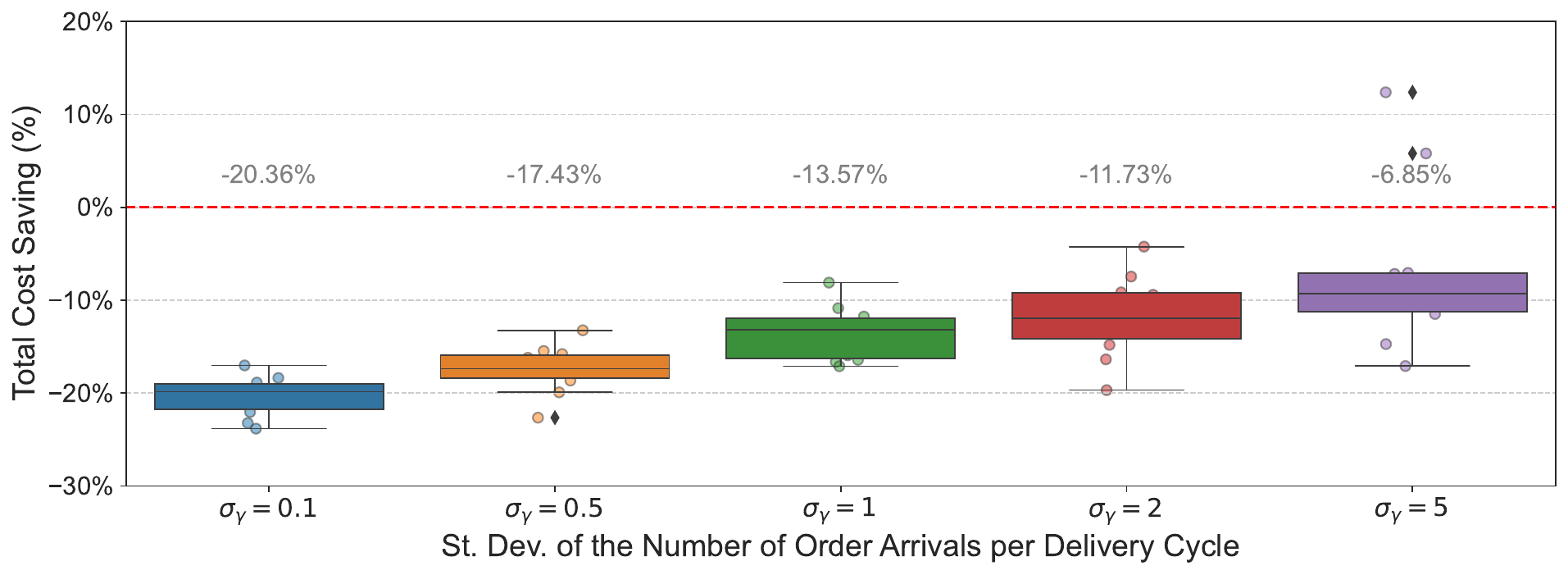}
    	\caption{Effect of uncertainties in order arrival intensity}
    	\label{sen_num_total}
        \end{subfigure}\\
        \begin{subfigure}[b]{.9\textwidth}
            \centering
            \includegraphics[width=1.0\textwidth]{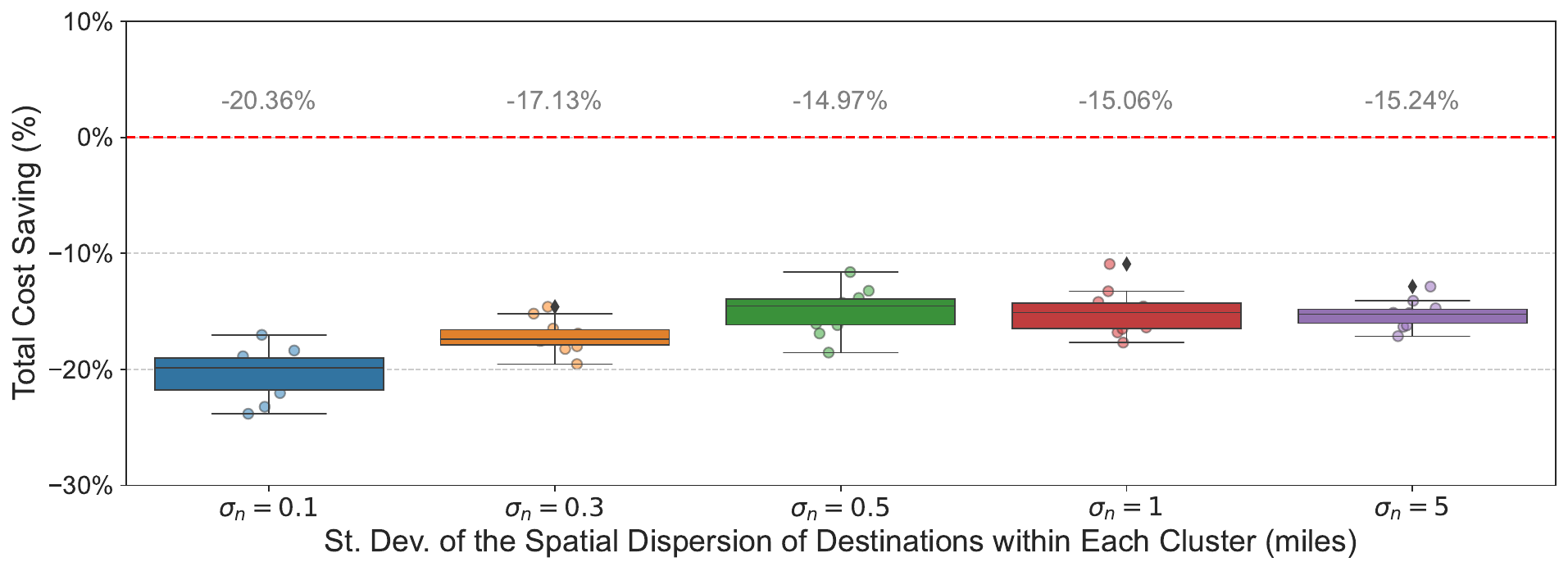}
    	\caption{Effect of uncertainties in destinations' spatial dispersion}
    	\label{sen_loc_total}
        \end{subfigure}\\
        \begin{subfigure}[b]{.9\textwidth}
            \centering
            \includegraphics[width=1.0\textwidth]{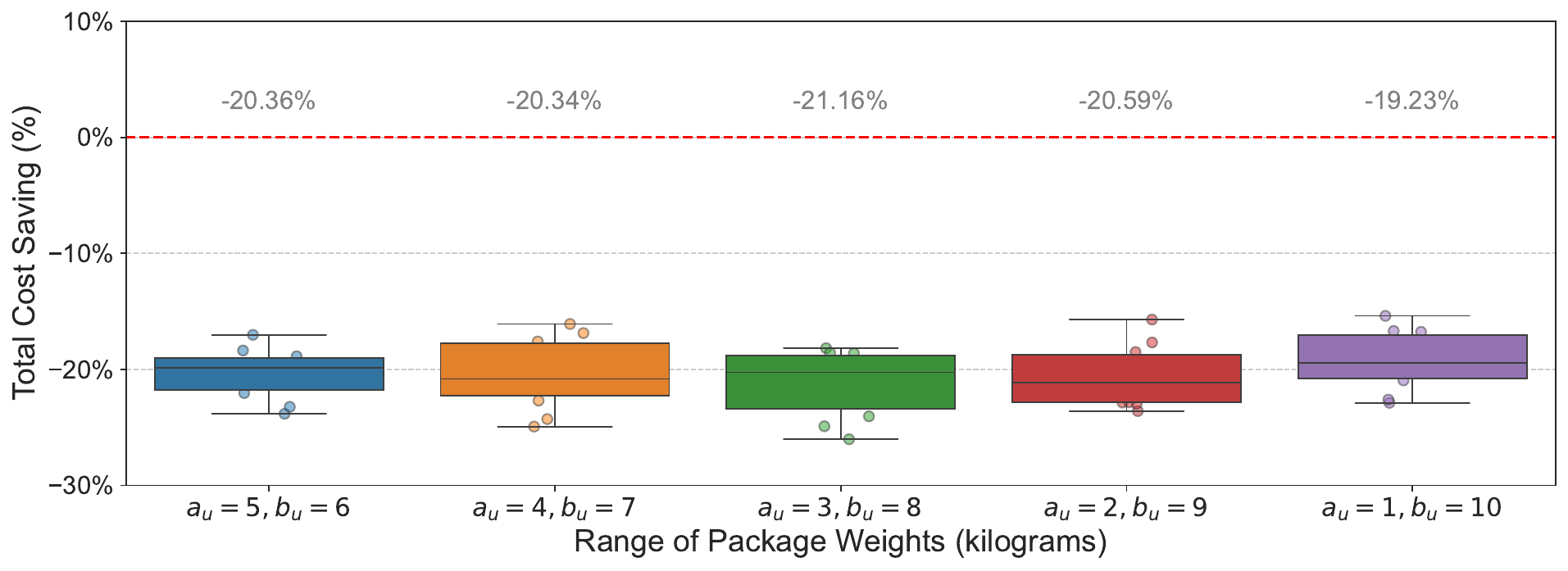}
    	\caption{Effect of uncertainties in package weights}
    	\label{sen_weight_total}
        \end{subfigure}
        \caption{Cost reductions achieved by the two-stage model compared to single-stage counterpart in response to different sources and levels of uncertainties (10 replications per scenario).}
        \label{fig_lunch peak_sensitivity}
\end{figure}

Building upon the base scenario, we then adjust the variances of various system inputs to evaluate the performance of the two-stage order-dispatch model under distinct aspects and levels of uncertainties. Specifically, our examination focuses on three aspects of uncertainties related to the order arrival intensity, their spatial dispersion, and package weights, respectively. \Cref{fig_lunch peak_sensitivity} displays the changes in total costs resulting from the two-stage versus single-stage models in response to these three sources of uncertainties.

\textit{Effect of uncertainties in order arrival intensity.}\quad Figure \ref{sen_num_total} shows that the total cost savings brought by the two-stage stochastic program over its single-stage counterpart diminish as the uncertainty increases in the intensity of delivery orders. As the standard deviation in the number of per-cycle order arrivals $\sigma_\gamma$ increases to 5 (exceeding the mean order arrival intensity), the cost-saving of the two-stage model decreases from over 20\% to 6\%. Occasionally, the two-stage model might produce solutions that are significantly suboptimal compared to those of the single-stage model. This highlights that when the demand intensity realizations are subject to high uncertainties (i.e., noise), the benefits of considering future information are diminished.

\textit{Effect of uncertainties in destinations' spatial dispersion.}\quad According to Figure \ref{sen_loc_total}, the cost savings achieved through the two-stage program initially narrow from 20\% to 15\% as the standard deviation in order destinations' spatial dispersion $\sigma_n$ increases from 0.1 to 0.5. Eventually, it stabilizes at around 15\% despite $\sigma_n$ increases further. This trend differs from the monotonic performance reductions in response to the uncertainties in the number of order arrivals. Indeed, increasing uncertainties in order locations lowers the performance of the two-stage program over the trip-saving mechanism. Nevertheless, it retains its ability to strategically delay and batch orders to trigger trip reductions constantly. As a result, even when uncertainties in order destinations become significantly high, the two-stage stochastic program continues to outperform the single-stage counterpart, maintaining the 15\% margin in cost savings.

\textit{Effect of uncertainties in package weights.}\quad Within the provided range of uncertainties, increasing the variability in package weights does not seem to significantly impact the performance of our proposed two-stage model. As shown in \Cref{sen_weight_total}, across all five levels of variability examined (by ranging $[a_u, b_u]$ from $[5, 6]$ to $[1, 10]$), the stochastic programming approach consistently achieves a 20\% cost reduction compared to the single-stage model. One possible reason is that the load capacity of the drones in our case can accommodate two or more orders per delivery trip, regardless of their variances in package weights. This provides a certain buffering for order batching, which helps to hedge against uncertainties in the fulfillment with different package weights.

\section{Conclusion}\label{sec:conclusion}

This paper investigated the dynamic order fulfillment problem in the context of last-mile drone delivery with demand uncertainty. A two-stage stochastic program was formulated according to the rolling horizons of the order fulfillment process, with the first stage considering the immediate delivery cycle and the second stage capturing several future cycles. Various sources of demand uncertainty arising from stochastic order arrivals, delivery destinations, and fulfillment urgency were incorporated. A branch-and-cut algorithm that incorporates two-tiered heuristics is devised based on the integer L-shaped method by \cite{laporte1993integer}, adapting the idea of an improved alternating strategy from \cite{angulo2016improving} and the machine-learning-based heuristics from \cite{larsen2023fast}, to effectively solve the two-stage model. The original exact L-shaped method suffers from great computational inefficiencies due to the vehicle routing decisions in our problem, making it unsuitable for the on-demand setting with real-time decision-making. Our enhanced L-shaped algorithm presented remarkable performance gains credited to: (i) the slimming of optimality cuts and the pruning of branching trees, reducing the required number of branching tree nodes; (ii) the substitution of exact solutions of the second-stage subproblem with fast heuristics for generating optimality cuts; (iii) the use of an alternating strategy that wisely adds `greedy heuristic' and `augmented heuristic' cuts. The implications of deploying the stochastic program for operationalizing on-demand drone delivery were examined through systematic numerical experiments.

Future research can further develop the current model in several meaningful ways. Firstly, while our current enhanced L-shaped algorithm shows promising improvements in computational efficiency, it still encounters challenges as the problem complexity grows exponentially in response to the number of pending orders in the first stage. One potential solution could involve `divide-and-conquer' strategies, such as partitioning the service area into smaller regions or pre-clustering orders by small neighborhoods. Secondly, concerning the operations of drones, the model could incorporate additional realistic protocols. A crucial aspect of drones is their battery-constrained limited delivery range compared to traditional gas-powered vehicles. Our current model assumes that each drone starts with a fully charged battery at the beginning of each cycle. An enhancement could involve integrating the charging process into the model, particularly for peak hours with heavy tasks and insufficient capacity for battery charging. Thirdly, the current work lays the foundation by considering a single-depot context, whilst future research could expand and generalize the model to multi-depot settings, which could be more applicable for larger spatial markets. For instance, exploring the incorporation of satellite hubs for dispatching, receiving, and charging drones could be relevant. Lastly, we could investigate the interactions between the logistical system discussed and the urban traffic system, taking into account the impacts from/on surface and aerial congestion.

\section*{Acknowledgement}
The work described in this paper was partly supported by the National Science Foundation (CMMI-2443338).

\section*{Declaration of generative Al and Al-assisted technologies in the writing process}
During the preparation of this work the authors used ChatGPT in order to enhance the readibility of the paper. After using this tool/service, the authors reviewed and edited the content as needed and take full responsibility for the content of the published article.

\newpage
\section*{Nomenclature}
\begin{table}[htb!]
        \footnotesize
	\centering
	\caption{Notation list of sets, variables, parameters and functions}
	\label{notatoin1}
	\begin{tabular}{c l} 
	\hline
	Notation & Description \\
	\hline\hline
	\multicolumn{2}{l}{\textbf{\emph{Sets}}}\\
	$m \in M$ & Existing requests in the current delivery cycle\\
	$n \in N$ & Potential requests arriving in the future cycles\\
        $(i,j) \in A (\tilde{A})$ & Feasible arcs for drones in the first stage (second stage) \\
	$k \in K$ & Fleet of drones\\
	$t \in \netT$ & Future delivery cycle inside the planning horizon\\	
	$\omega \in \Omega$ & Scenario for potential realizations of the future cycles\\
	\hline
	
	\multicolumn{2}{l}{\textbf{\emph{Common Parameters}}}\\
	$c_m^t$ & Penalty for a existing task $m$ delivered in the $t^{\text{th}}$ cycle\\
	$d_{ij}$ & Delivery distance between nodes $i$ and $j$\\
        $q_m$ & Package weight of request $m$ \\
        $Q$ & Load capacity of individual drone \\
        $E$ & Battery capacity of individual drone \\
        $\mathcal{M}$ & Constant of extremely large value\\
        $\alpha, \beta, \gamma, \kappa$ & Coefficients for different cost items \\
	
	\multicolumn{2}{l}{\textbf{\emph{Parameters for Scenario $\xi$}}}\\
	$c^{t,\xi}_n$ & Penalty for a second-stage task $n$ delivered in the $t^{\text{th}}$ cycle\\
        $q^{\xi}_n$ & Package weight of request $n$ \\
	$a^{\xi}_n$ & First decision epoch for request $n$ following its arrival\\
	\hline
	
	\multicolumn{2}{l}{\textbf{\emph{First-Stage Variables}}}\\
	$y_m^0$ & Whether existing order $m$ is fulfilled in the current delivery cycle\\
	$x_{ij}^{0,k}$ & Whether drone $k$ travels from node $i$ to $j$ in the current delivery cycle\\
        $r^{0,k}$ & Whether drone $k$ is dispatched in the current delivery cycle \\
        $u_i^{0,k}$ & Carrying load of drone $k$ after visiting node $i$, kilograms  \\
        $u_O^{0,k}, u_D^{0,k}$ & Starting and ending loads of drone $k$ in the current delivery cycle, kilograms \\
        $v_i^{0,k}$ & State of charge of drone $k$ after visiting node $i$, kilowatt-hours\\
        $v_O^{0,k}, v_D^{0,k}$ & Starting and ending states of charge for drone $k$ in the current delivery cycle, kilowatt-hours\\

	\multicolumn{2}{l}{\textbf{\emph{Second-Stage Variables for Scenario $\xi$}}}\\
        $y_m^{t,\xi} \ (y_n^{t,\xi})$ & Whether existing order $m$ (future order $n$) would be fulfilled in the $t^{\text{th}}$ delivery cycle\\
        $x_{ij}^{t,k,\xi}$ & Whether drone $k$ would travel from node $i$ to $j$ in the $t^{\text{th}}$ delivery cycle\\
        $r^{t,k,\xi}$ & Whether drone $k$ is dispatched in the $t^{\text{th}}$ delivery cycle \\
        $u_i^{t,k,\xi}$ & Carrying load of drone $k$ after visiting node $i$ in the $t^{\text{th}}$ delivery cycle, kilograms\\
        $u_O^{t,k,\xi}, u_D^{t,k,\xi}$ & Starting and ending loads of drone $k$ in the $t^{\text{th}}$ delivery cycle, kilograms\\
        $v_i^{t,k,\xi}$ & State of charge of drone $k$ after visiting node $i$ in the $t^{\text{th}}$ delivery cycle, kilowatt-hours\\
        $v_O^{t,k,\xi}, v_D^{t,k,\xi}$ & Starting and ending states of charge for drone $k$ in the $t^{\text{th}}$ delivery cycle, kilowatt-hours\\
	$z_n^{\xi}$ & Whether request $n$ would be handled by the virtual drone upon generation\\
	\hline
	
	\multicolumn{2}{l}{\textbf{\emph{Functions}}}\\
        $P(\cdot)$ & Energy consumption function of individual drone based on its carrying load\\
        $\mathbb{Q}(\cdot, \xi)$ & Second-stage value function given scenario $\xi$\\
        $\mathbb{L}(\cdot)$ & Second-stage expected value (recourse) function \\
	\hline
	\end{tabular}
\end{table}
\clearpage

\begin{appendices}
    \section{Proof of Propositions}\label{sec:appd_proof}
\setcounter{equation}{0}
\setcounter{figure}{0}
\setcounter{table}{0}
\renewcommand\theequation{A\arabic{equation}}
\renewcommand\thefigure{A\arabic{figure}}
\renewcommand\thetable{A\arabic{table}}

This appendix presents the proofs for the propositions introduced in \Cref{sec:int_L-shaped_method}.

\subsection{Proof of Proposition \ref{proposition_1}}

\begin{proof}
We prove that \Cref{opt_cut_1} for $(\bm{y^{0}}, \theta)$ constitutes a valid optimality cut that would not exclude the optimal solution $(\bm{y^{0,*}}, \mathbb{L}(\bm{y^{0,*}}))$. For any arbitrary iteration $\eta$, if $\bm{y^{0,\eta}}$ achieves the optimal solution $\bm{y^{0,*}}$, we have
\begin{align*}
    \left(\sum_{m \in S^\eta} y_m^{0,*} - \sum_{m \notin S^\eta} y_m^{0,*} \right) = \left(\sum_{m \in S^\eta} y_m^{0,\eta} - \sum_{m \notin S^\eta} y_m^{0,\eta} \right) = \left| S^\eta \right|,
\end{align*}
with which the optimality cut at $\bm{y^{0,*}}$ reduces to 
\begin{align*}
    \theta &\geq (\mathbb{L}(\bm{y^{0,\eta}})-L) \left(\sum_{m \in S^\eta} y_m^{0,*} - \sum_{m \notin S^\eta} y_m^{0,*} \right) - (\mathbb{L}(\bm{y^{0,\eta}})-L)(\left| S^\eta \right| -1) + L \\
    & = (\mathbb{L}(\bm{y^{0,\eta}})-L) \left| S^\eta \right| - (\mathbb{L}(\bm{y^{0,\eta}})-L)(\left| S^\eta \right| -1) + L \\
    &= \mathbb{L}(\bm{y^{0,\eta}}).
\end{align*}
Clearly, the solution $\theta = \mathbb{L}(\bm{y^{0,*}})$ is tightly included in $\theta \geq \mathbb{L}(\bm{y^{0,\eta}})$.

If $\bm{y^{0,\eta}}$ differs from the optimal solution $\bm{y^{0,*}}$, we have $\left(\sum_{m \in S^\eta} y_m^{0,*} - \sum_{m \notin S^\eta} y_m^{0,*} \right)$ being strictly less than $\left| S^\eta \right|$. As all elements in $\bm{y^{0,*}}$ are integral, the right-hand side (RHS) of \Cref{opt_cut_1} at the optimal solution $\bm{y^{0,*}}$ can be expressed as
\begin{align*}
    \text{RHS} &= (\mathbb{L}(\bm{y^{0,\eta}})-L) \left(\sum_{m \in S^\eta} y_m^{0,*} - \sum_{m \notin S^\eta} y_m^{0,*} \right) - (\mathbb{L}(\bm{y^{0,\eta}})-L)(\left| S^\eta \right| -1) + L \\
    & \leq (\mathbb{L}(\bm{y^{0,\eta}})-L) (\left| S^\eta \right| - 1) - (\mathbb{L}(\bm{y^{0,\eta}})-L)(\left| S^\eta \right| -1) + L \\
    &= L.
\end{align*}
Since $\mathbb{L}(\bm{y^{0,*}})\geq L$, the solution $(\bm{y^{0,*}}, \mathbb{L}(\bm{y^{0,*}}))$ is feasible to the optimality cut. In this case, the cut is valid but likely not tight.

\end{proof}

\subsection{Proof of Proposition \ref{proposition_2}}

\begin{proof}
Because all the cost terms in the second-stage objective in Eq. (\ref{opt_model_1}b) are non-negative, a naive lower bound would be $L=0$. Here, we establish a more reasonable (and useful) lower bound $\mathbb{L}(\bm{y^{0}}=\bm{1})$ for $L$ by setting the first-stage variables $\bm{y^0}$ to $1$, suggesting all existing orders are fulfilled in the current delivery cycle with no remaining orders passed to the future. We note that $\bm{y^{0}}=\bm{1}$ may be infeasible for the first stage problem. This proposition, however, is to prove that such a strategy produces a valid lower bound for the second-stage recourse value, regardless of whether this bound is attainable or not.

For a given scenario $\omega$, consider an arbitrary order fulfillment plan $\{\bm{y^0, y_M^{t,\omega}}, \netZ^{t,\omega}\}_{t\in\netT\cup \{T+1\}}$ that is feasible for the second-stage problem and yields the second-stage objective value $\hat{\mathbb{Q}}(\bm{y^0}, \bm{y_M^{t, \omega}}, \bm{\netZ^{t, \omega}})$. In this plan, $\bm{y^0}$ and $\bm{y_M^{t, \omega}}$ indicate the fulfillment schedule for existing orders, with $y_m^0 + \sum_{t\in \netT \cup \{T+1\}} y_m^{t, \omega} = 1$ for all $m\in M$. Besides, $\netZ^{t, \omega}$ includes all the decision variables in the second stage, other than $\bm{y_M^{t, \omega}}$, defining the fulfillment schedules of future orders as well as the routing of drones in future delivery cycles in scenario $\omega$, i.e., $\netZ^{t, \omega} = \{\bm{y_N^{t, \omega}}, \bm{z^{t, \omega}}, \bm{x^{t, \omega}}, \bm{r^{t, \omega}}\}_{t\in\netT\cup \{T+1\}}$.

Suppose that according to the plan $\{\bm{y^0, y_M^{t,\omega}}, \netZ^{t,\omega}\}_{t\in\netT\cup \{T+1\}}$, some existing orders are deferred to future delivery cycles, i.e. $\exists m \in M$, $y_m^0 = 0$, and $y_m^{t,\omega} = 1$ for certain $t\in \netT\cup \{T+1\}$. Then, we consider an alternative plan $\{\bm{1, 0}, \netZ^{t,\omega}\}_{t\in\netT\cup \{T+1\}}$ that fulfills all existing orders in the current cycle without deferral. This alternative plan keeps the fulfillment schedules of all future orders unchanged and maintains the same delivery routes for drones as the original plan. With the fulfillment of existing orders moved up, the carrying load of drones in future delivery cycles either remains the same or decreases, ensuring that the original routing remains feasible under the alternative plan. Consequently, all cost terms related to order delay penalties for second-stage orders $N$ and fixed costs of drone dispatches remain unchanged, while the order delay penalties for existing orders $M$ and the energy consumption of drones are reduced for executing the alternative plan. Thus, we have $\hat{\mathbb{Q}}(\bm{1, 0}, \netZ^{t, \omega}) \leq \hat{\mathbb{Q}}(\bm{y^0}, \bm{y_M^{t, \omega}}, \bm{\netZ^{t, \omega}})$ for any second-stage order fulfillment plan $\{\bm{y^0, y_M^{t,\omega}}, \netZ^{t,\omega}\}_{t\in\netT\cup \{T+1\}}$. Then, with the front-fulfillment plan for existing orders in place, resolving the optimal order fulfillment and drone routing plan $\bar{\netZ}^{t, \omega, *}$ gives rise to $\hat{\mathbb{Q}}(\bm{1, 0}, \bar{\netZ}^{t, \omega, *}) \leq \hat{\mathbb{Q}}(\bm{1, 0}, \netZ^{t, \omega})$.

The above inequalities hold for any given scenario $\omega$. Taking the expectations across all possible scenarios further leads to the following relationships,
\begin{align*}
    \frac{1}{|\Omega|}\sum_{\omega \in \Omega}  \mathbb{Q}(\bm{y^0},\omega) &\geq \frac{1}{|\Omega|}\sum_{\omega \in \Omega}  \hat{\mathbb{Q}}(\bm{1, 0}, \netZ^{t, \omega}) \\
    & \geq \frac{1}{|\Omega|}\sum_{\omega \in \Omega}  \hat{\mathbb{Q}}(\bm{1, 0}, \netZ^{t, \omega,*}) \\
    & = \frac{1}{|\Omega|}\sum_{\omega \in \Omega} \mathbb{Q}(\bm{1},\omega) \\
    & = \mathbb{L}(\bm{1}).
\end{align*}
This concludes that $\mathbb{L}(\bm{y^{0}}=\bm{1})$ serves as a valid lower bound on the second-stage recourse costs for all possible $\bm{y^{0}}$. 

\end{proof}

\subsection{Proof of Proposition \ref{proposition_3}}

\begin{proof}
Suppose there exists a hyperparameter $\nu$ such that $\nu \mathbb{L}^{gred}(\bm{y^{0}}) < \mathbb{L}(\bm{y^{0}})$. Then, for a particular iteration $\eta$, if the condition $\theta^{\eta} < \nu \mathbb{L}^{gred}(\bm{y^{0,\eta}})$ is met, the obtained $\bm{y^{0,\eta}}$ is not optimal, as the optimal solution must satisfy $\theta^{\eta} \geq \mathbb{L}(\bm{y^{0,\eta}})$ \citep{laporte1993integer}. As we shown in the proof of Proposition \ref{proposition_1}, a non-optimal but feasible $\bm{y^{0,\eta}}$ satisfies $\left(\sum_{m \in S^\eta} y_m^{0,\eta} - \sum_{m \notin S^\eta} y_m^{0,\eta} \right) \leq \left| S^\eta \right| - 1$.

Next, given $\bm{y^{0,\eta}}$, for the second-stage problem obtained using `augmented heuristic' and `greedy heuristic' methods, we have the inequality $\mathbb{L}(\bm{y^{0,\eta}}) \leq \mathbb{L}^{aug}(\bm{y^{0,\eta}}) \leq \mathbb{L}^{gred}(\bm{y^{0,\eta}})$. Consequently, between the derived `augmented cut' and `greedy cut', their right-hand side terms as a function of $y_m^0$ are related as follows:
\begin{align*}
     & \text{RHS}^{gred} (y_m^0) - \text{RHS}^{aug} (y_m^0) \\
     =& \left(\mathbb{L}^{gred}(\bm{y^{0,\eta}})-\mathbb{L}^{aug}(\bm{y^{0,\eta}})\right) \left(\sum_{m \in S^\eta} y_m^0 - \sum_{m \notin S^\eta} y_m^0 \right) - \left(\mathbb{L}^{gred}(\bm{y^{0,\eta}})-\mathbb{L}^{aug}(\bm{y^{0,\eta}})\right)(\left| S^\eta \right| -1) \\
     \leq& \left(\mathbb{L}^{gred}(\bm{y^{0,\eta}})-\mathbb{L}^{aug}(\bm{y^{0,\eta}})\right)(\left| S^\eta \right| -1) - \left(\mathbb{L}^{gred}(\bm{y^{0,\eta}})-\mathbb{L}^{aug}(\bm{y^{0,\eta}})\right)(\left| S^\eta \right| -1) \\
     =& 0.
\end{align*}
The `augmented cut'--- $\theta \geq \text{RHS}^{aug}(y_m^0)$ ---is thus tighter than the `greedy cut'--- $\theta \geq \text{RHS}^{gred}(y_m^0)$ ---, and adding the latter would not exclude the solutions yielded by imposing the former. 
\end{proof}

    \section{Heuristic Methods for Solving the Subproblem (\textbf{SP}) with Vehicle Routing Constraints}\label{sec:appd_subproblem_acc}
\setcounter{equation}{0}
\setcounter{figure}{0}
\setcounter{table}{0}
\renewcommand\theequation{B\arabic{equation}}
\renewcommand\thefigure{B\arabic{figure}}
\renewcommand\thetable{B\arabic{table}}

This appendix details the heuristic methods employed to solve (\textbf{SP}) using the Google OR-Tools' vehicle routing package. For each scenario $\omega \in \Omega$, upon the first-stage solution $\bm{y^{0,\eta}}$, we unfold the subproblem (\textbf{SP}) as a capacitated VRP with time-window constraints. This is achieved by aggregating all pending and future orders and incorporating penalties for dropping visits. 

We discretize time windows based on the $T$ delivery cycles, and have each vehicle operate only within one specific time window. Thus, a single vehicle in service throughout the $T$ cycles is copied $T$ times, resulting in a total of $|K|\times T$ vehicles for the Google OR-Tool's API. Moreover, in the dynamic fulfillment system, orders fall into three categories: unfulfilled orders inherited from over $T$ cycles ago, unfulfilled orders from previous cycles but less than $T$ cycles, and potential orders in the future. Let $\tau(t)$ denote the time of transition between the $t^{\text{th}}$ and $(t+1)^{\text{th}}$ cycles, representing the end time of the $t^{\text{th}}$ cycle and the start time of the $(t+1)^{\text{th}}$ cycle. Then, orders in the three categories are assigned with soft time windows of $[0,\tau(1))$, $[0, \tau(T-a_m))$, and $[\tau(a_m), \tau{T})$, respectively, where $a_m$ denotes the arrival cycle of a particular order $m$. Unlike the strict time windows for vehicles, these soft time windows are defined to associate penalties, $c_m^{T+1}$, with fulfillment delays. Given the capacity provided by the Google OR-Tool, the objective function submitted to the API accounts for the penalties arising from delayed orders, vehicle travel distances, and the frequencies of vehicle dispatch. Energy consumption is treated separately and will be addressed later. Two advanced solution-searching strategies, namely ``PATH\_CHEAPEST\_ARC'' and ``GUIDED\_LOCAL\_SEARCH'', are utilized as an enhancement to the constraint-programming-based VRP solver. The former strategy swiftly generates a feasible solution using a greedy approach, while the latter is designed to escape from local minima and continues searching for improvements\footnote{\url{https://developers.google.com/optimization/routing/routing_options}}. In this study, the `greedy heuristic' solely implements ``PATH\_CHEAPEST\_ARC'', while the `augmented heuristic' combines both ``PATH\_CHEAPEST\_ARC'' and ``GUIDED\_LOCAL\_SEARCH''.

As mentioned, the energy consumption aspect in (\textbf{SP}) is not addressed by the Google OR-Tool. Its primary impact on the discussed decisions relates to the sequencing of package deliveries, reflecting a trade-off between selecting orders with heavier packages and those destined for locations at a distance. Accordingly, we propose a fix to adjust solutions based on the metric of energy consumption. First, we receive the API's output of order assignments that can be feasibly serviced by each vehicle at each delivery cycle. Then, on top of this resulting vehicle-order assignment, we proceed to solve a travelling salesman problem for each individual vehicle to fine-tune their routing, particularly considering energy consumption. Such a remedy, involving the separate treatment for cost factors, may be considered somewhat greedy. Nevertheless, as reported in \Cref{tab_sub_comparison_4,tab_sub_comparison_6,tab_sub_comparison_8}, this approach significantly reduces the computational time with only minor compromise in solution quality.

\begin{table}[b!]
    \small
    \centering
    \caption{Comparisons between heuristic and exact methods for solving (\textbf{SP}) \\ ($\left | M \right |=4$, $\left | N \right |=4$, $T=2$, $\left | K \right |=2$, $\left |\Omega \right |=10$)}\label{tab_sub_comparison_4}
\begin{tabular}{c|ccc|ccc}
\hline
\multirow{2}{*}{Instance} & \multicolumn{3}{c|}{Gap on Minimum Recourse Costs $\mathbb{L}(\bm{y^{0,\eta}}$)}                                          & \multicolumn{3}{c}{Running Time (seconds)}                                        \\ \cline{2-7} 
                      & \multicolumn{1}{c|}{Augmented Heu.} & \multicolumn{1}{c|}{Greedy Heu.} & Gurobi 9.5 & \multicolumn{1}{c|}{Augmented Heu.} & \multicolumn{1}{c|}{Greedy Heu.} & Gurobi 9.5 \\ \hline
1                     & \multicolumn{1}{c|}{0.46\%}    & \multicolumn{1}{c|}{0.46\%}         & 0       & \multicolumn{1}{c|}{1.02}      & \multicolumn{1}{c|}{0.02}           & 2.19       \\
2                     & \multicolumn{1}{c|}{1.15\%}    & \multicolumn{1}{c|}{1.55\%}         & 0       & \multicolumn{1}{c|}{1.01}      & \multicolumn{1}{c|}{0.01}           & 6.66       \\
3                     & \multicolumn{1}{c|}{0.07\%}    & \multicolumn{1}{c|}{2.87\%}         & 0       & \multicolumn{1}{c|}{1.01}      & \multicolumn{1}{c|}{0.01}           & 0.10       \\
4                     & \multicolumn{1}{c|}{0.93\%}    & \multicolumn{1}{c|}{2.86\%}         & 0       & \multicolumn{1}{c|}{1.01}      & \multicolumn{1}{c|}{0.01}           & 3.97       \\
5                     & \multicolumn{1}{c|}{0.55\%}    & \multicolumn{1}{c|}{1.77\%}         & 0       & \multicolumn{1}{c|}{1.01}      & \multicolumn{1}{c|}{0.02}           & 5.12       \\
6                     & \multicolumn{1}{c|}{2.40\%}    & \multicolumn{1}{c|}{2.51\%}         & 0       & \multicolumn{1}{c|}{1.01}      & \multicolumn{1}{c|}{0.01}           & 0.20       \\
7                     & \multicolumn{1}{c|}{0.45\%}    & \multicolumn{1}{c|}{0.45\%}         & 0       & \multicolumn{1}{c|}{1.01}      & \multicolumn{1}{c|}{0.01}           & 0.29       \\
8                     & \multicolumn{1}{c|}{0.17\%}    & \multicolumn{1}{c|}{0.17\%}         & 0       & \multicolumn{1}{c|}{1.01}      & \multicolumn{1}{c|}{0.01}           & 1.25       \\
9                     & \multicolumn{1}{c|}{0.18\%}    & \multicolumn{1}{c|}{1.73\%}         & 0       & \multicolumn{1}{c|}{1.01}      & \multicolumn{1}{c|}{0.01}           & 0.46       \\
10                    & \multicolumn{1}{c|}{0.00\%}    & \multicolumn{1}{c|}{0.00\%}         & 0       & \multicolumn{1}{c|}{1.01}      & \multicolumn{1}{c|}{0.01}           & 0.20       \\ \hline
Avg.               & \multicolumn{1}{c|}{0.64\%}    & \multicolumn{1}{c|}{1.44\%}         & 0       & \multicolumn{1}{c|}{1.01}      & \multicolumn{1}{c|}{0.01}           & 2.04       \\
St. Dev.               & \multicolumn{1}{c|}{0.72\%}    & \multicolumn{1}{c|}{1.11\%}         & 0       & \multicolumn{1}{c|}{0.00}      & \multicolumn{1}{c|}{0.00}           & 2.39       \\
Max.                   & \multicolumn{1}{c|}{2.40\%}    & \multicolumn{1}{c|}{2.87\%}         & 0       & \multicolumn{1}{c|}{1.02}      & \multicolumn{1}{c|}{0.02}           & 6.66       \\ \hline
\end{tabular}
\end{table}

A series of experiments are conducted to evaluate the performances of the `greedy heuristic' and the `augmented heuristic' methods in solving the subproblem. The solution quality of these methods is measured by the gap on minimum recourse costs, which compares the obtained objective values against the exact solutions provided by Gurobi, whilst efficiency is measured by the running time for executing instances. \Cref{tab_sub_comparison_4,tab_sub_comparison_6,tab_sub_comparison_8} present the results of 30 instances corresponding to three different parametric setups. Each instance computing $\left |\Omega \right |=10$ scenarios and reports the average recourse costs obtained and the average running time per scenario. The results clearly demonstrate the significant superiority of heuristic methods in problem-solving. Even for relatively small subproblems, such as a VRP with 12 nodes (see \Cref{tab_sub_comparison_6} with $\left | M \right |=6$ and $\left | N \right |=6$), Gurobi takes as long as 189 seconds to derive the exact solution. The challenge escalates with larger problems; for instance, with a 16-node problem (see \Cref{tab_sub_comparison_8} with $\left | M \right |=8$ and $\left | N \right |=8$), half of the instances require more than 10 minutes. In a computation involving 10 scenarios, this translates to an execution time of 100 minutes to generate one optimality cut without parallel computing, a timeline impractical for online problem-solving. In stark contrast, the `augmented heuristic' and `greedy heuristic' methods demonstrate remarkable efficacy, solving the same problems within 1 second and 0.03 seconds, respectively, while maintaining satisfactory optimality gaps.

\begin{table}[!ht]
    \small
    \centering
    \caption{Comparisons between heuristic and exact methods for solving (\textbf{SP}) \\ ($\left | M \right |=6$, $\left | N \right |=6$, $T=2$, $\left | K \right |=2$, $\left |\Omega \right |=10$)}\label{tab_sub_comparison_6}
\begin{tabular}{c|ccc|ccc}
\hline
\multirow{2}{*}{Instance} & \multicolumn{3}{c|}{Gap on Minimum Recourse Costs $\mathbb{L}(\bm{y^{0,\eta}}$)}                                      & \multicolumn{3}{c}{Running Time (seconds)}                                        \\ \cline{2-7} 
                      & \multicolumn{1}{c|}{Augmented Heu.} & \multicolumn{1}{c|}{Greedy Heu.} & Gurobi 9.5 & \multicolumn{1}{c|}{Augmented Heu.} & \multicolumn{1}{c|}{Greedy Heu.} & Gurobi 9.5 \\ \hline
1                     & \multicolumn{1}{c|}{1.36\%}    & \multicolumn{1}{c|}{2.00\%}         & 0          & \multicolumn{1}{c|}{1.02}      & \multicolumn{1}{c|}{0.03}           & 333.87     \\
2                     & \multicolumn{1}{c|}{0.09\%}    & \multicolumn{1}{c|}{0.16\%}         & 0          & \multicolumn{1}{c|}{1.02}      & \multicolumn{1}{c|}{0.02}           & 108.21     \\
3                     & \multicolumn{1}{c|}{0.81\%}    & \multicolumn{1}{c|}{1.54\%}         & 0          & \multicolumn{1}{c|}{1.02}      & \multicolumn{1}{c|}{0.03}           & 307.61     \\
4                     & \multicolumn{1}{c|}{0.05\%}    & \multicolumn{1}{c|}{0.46\%}         & 0          & \multicolumn{1}{c|}{1.02}      & \multicolumn{1}{c|}{0.02}           & 117.82     \\
5                     & \multicolumn{1}{c|}{0.09\%}    & \multicolumn{1}{c|}{4.56\%}         & 0          & \multicolumn{1}{c|}{1.02}      & \multicolumn{1}{c|}{0.02}           & 112.08     \\
6                     & \multicolumn{1}{c|}{0.31\%}    & \multicolumn{1}{c|}{0.31\%}         & 0          & \multicolumn{1}{c|}{1.02}      & \multicolumn{1}{c|}{0.02}           & 21.16      \\
7                     & \multicolumn{1}{c|}{0.33\%}    & \multicolumn{1}{c|}{0.77\%}         & 0          & \multicolumn{1}{c|}{1.01}      & \multicolumn{1}{c|}{0.02}           & 2.66       \\
8                     & \multicolumn{1}{c|}{0.45\%}    & \multicolumn{1}{c|}{2.63\%}         & 0          & \multicolumn{1}{c|}{1.02}      & \multicolumn{1}{c|}{0.02}           & 379.52     \\
9                     & \multicolumn{1}{c|}{0.54\%}    & \multicolumn{1}{c|}{2.20\%}         & 0          & \multicolumn{1}{c|}{1.02}      & \multicolumn{1}{c|}{0.02}           & 12.61      \\
10                    & \multicolumn{1}{c|}{1.65\%}    & \multicolumn{1}{c|}{4.41\%}         & 0          & \multicolumn{1}{c|}{1.02}      & \multicolumn{1}{c|}{0.03}           & 499.27     \\ \hline
Avg.               & \multicolumn{1}{c|}{0.57\%}    & \multicolumn{1}{c|}{1.90\%}         & 0          & \multicolumn{1}{c|}{1.02}      & \multicolumn{1}{c|}{0.02}           & 189.48     \\
St. Dev.               & \multicolumn{1}{c|}{0.55\%}    & \multicolumn{1}{c|}{1.60\%}         & 0          & \multicolumn{1}{c|}{0.00}      & \multicolumn{1}{c|}{0.00}           & 176.11     \\
Max.                   & \multicolumn{1}{c|}{1.65\%}    & \multicolumn{1}{c|}{4.56\%}         & 0          & \multicolumn{1}{c|}{1.02}      & \multicolumn{1}{c|}{0.03}           & 499.27     \\ \hline
\end{tabular}
\end{table}

\begin{table}[!ht]
    \small
    \centering
    \caption{Comparisons between heuristic and exact methods for solving (\textbf{SP}) \\ ($\left | M \right |=8$, $\left | N \right |=8$, $T=2$, $\left | K \right |=2$, $\left |\Omega \right |=10$)}\label{tab_sub_comparison_8}
\begin{tabular}{c|ccc|ccc}
\hline
\multirow{2}{*}{Instance} & \multicolumn{3}{c|}{Gap on Minimum Recourse Costs $\mathbb{L}(\bm{y^{0,\eta}}$)}                                          & \multicolumn{3}{c}{Running Time (seconds)}                                 \\ \cline{2-7} 
                      & \multicolumn{1}{c|}{Augmented Heu.} & \multicolumn{1}{c|}{Greedy Heu.} & Gurobi 9.5 & \multicolumn{1}{c|}{Augmented Heu.} & \multicolumn{1}{c|}{Greedy Heu.} & Gurobi 9.5 \\ \hline
1                     & \multicolumn{1}{c|}{0.07\%}          & \multicolumn{1}{c|}{6.34\%}           & 0          & \multicolumn{1}{c|}{1.03}            & \multicolumn{1}{c|}{0.03}             & $>600$     \\
2                     & \multicolumn{1}{c|}{0.56\%}          & \multicolumn{1}{c|}{5.33\%}           & 0          & \multicolumn{1}{c|}{1.03}            & \multicolumn{1}{c|}{0.03}             & 549.02     \\
3                     & \multicolumn{1}{c|}{0.54\%}          & \multicolumn{1}{c|}{3.03\%}           & 0          & \multicolumn{1}{c|}{1.03}            & \multicolumn{1}{c|}{0.03}             & 587.34     \\
4                     & \multicolumn{1}{c|}{0.14\%}          & \multicolumn{1}{c|}{8.77\%}           & 0          & \multicolumn{1}{c|}{1.03}            & \multicolumn{1}{c|}{0.04}             & 583.47     \\
5                     & \multicolumn{1}{c|}{0.06\%}          & \multicolumn{1}{c|}{0.79\%}           & 0          & \multicolumn{1}{c|}{1.02}            & \multicolumn{1}{c|}{0.03}             & $>600$     \\
6                     & \multicolumn{1}{c|}{0.05\%}          & \multicolumn{1}{c|}{7.44\%}           & 0          & \multicolumn{1}{c|}{1.02}            & \multicolumn{1}{c|}{0.03}             & 462.84     \\
7                     & \multicolumn{1}{c|}{0.45\%}          & \multicolumn{1}{c|}{3.01\%}           & 0          & \multicolumn{1}{c|}{1.02}            & \multicolumn{1}{c|}{0.02}             & 145.08     \\
8                     & \multicolumn{1}{c|}{1.83\%}          & \multicolumn{1}{c|}{5.64\%}           & 0          & \multicolumn{1}{c|}{1.03}            & \multicolumn{1}{c|}{0.03}             & $>600$     \\
9                     & \multicolumn{1}{c|}{0.69\%}          & \multicolumn{1}{c|}{8.12\%}           & 0          & \multicolumn{1}{c|}{1.03}            & \multicolumn{1}{c|}{0.03}             & $>600$      \\
10                    & \multicolumn{1}{c|}{0.18\%}          & \multicolumn{1}{c|}{1.13\%}           & 0          & \multicolumn{1}{c|}{1.03}            & \multicolumn{1}{c|}{0.03}             & $>600$     \\ \hline
Avg.               & \multicolumn{1}{c|}{0.46\%}          & \multicolumn{1}{c|}{4.96\%}           & 0          & \multicolumn{1}{c|}{1.03}            & \multicolumn{1}{c|}{0.03}             & 533.21     \\
St. Dev.               & \multicolumn{1}{c|}{0.54\%}          & \multicolumn{1}{c|}{2.84\%}           & 0          & \multicolumn{1}{c|}{0.00}               & \multicolumn{1}{c|}{0.00}                & 143.06     \\
Max.                   & \multicolumn{1}{c|}{1.83\%}          & \multicolumn{1}{c|}{8.77\%}           & 0          & \multicolumn{1}{c|}{1.03}            & \multicolumn{1}{c|}{0.04}             & $>600$     \\ \hline
\end{tabular}
\end{table}
	\clearpage
\end{appendices}

\bibliographystyle{apalike} 
\bibliography{mybib} 



\end{document}